\documentclass[12pt]{article}
\usepackage{subeqn}
\usepackage{graphicx}
\usepackage{ifpdf}
\ifpdf  \fi
\usepackage{amsfonts,amssymb}
\usepackage{mathptmx,helvet,courier,makeidx,multicol,footmisc}
\usepackage[numbers]{natbib}
\bibpunct{(}{)}{;}{a}{,}{,}
\usepackage[bookmarksnumbered=true, pdfauthor={Wen-Long Jin}]{hyperref}

\newcommand{\commentout}[1]{}

\newcommand{\ba}{\begin{array}}
        \newcommand{\ea}{\end{array}}
\newcommand{\bc}{\begin{center}}
        \newcommand{\ec}{\end{center}}
\newcommand{\bdm}{\begin{displaymath}}
        \newcommand{\edm}{\end{displaymath}}
\newcommand{\bds} {\begin{description}}
        \newcommand{\eds} {\end{description}}%17Apr01
\newcommand{\ben}{\begin{enumerate}}
        \newcommand{\een}{\end{enumerate}}
\newcommand{\beq}{\begin{equation}}
        \newcommand{\eeq}{\end{equation}}
\newcommand{\bfg} {\begin{figure}[h]}
        \newcommand{\efg} {\end{figure}}%Nov 5,99
\newcommand{\bi} {\begin {itemize}}
        \newcommand{\ei} {\end {itemize}}
\newcommand{\bqn}{\begin{eqnarray}}
        \newcommand{\eqn}{\end{eqnarray}}
\newcommand{\bqs}{\begin{eqnarray*}}
        \newcommand{\eqs}{\end{eqnarray*}}
\newcommand{\bsl} {\begin{slide}[8.8in,6.7in]}
        \newcommand{\esl} {\end{slide}}
\newcommand{\bsq}{\begin{subequations}}
        \newcommand{\esq}{\end{subequations}}       
\newcommand{\bss} {\begin{slide*}[9.3in,6.7in]}
        \newcommand{\ess} {\end{slide*}}
\newcommand{\btb} {\begin {table}}
        \newcommand{\etb} {\end {table}}%Nov 10,99
\newcommand{\bvb}{\begin{verbatim}}
        \newcommand{\evb}{\end{verbatim}}

\newcommand{\m}{\mbox}
\newcommand {\der}[2] {{\frac {\m {d} {#1}} {\m{d} {#2}}}}

\newcommand {\pd}[2] {{\frac {\partial {#1}} {\partial {#2}}}}

\newcommand{\cas}[1]{{{\left \{ \ba #1 \ea \right. }}}

\newcommand{\reff}[1] {{{Figure \ref {#1}}}}
\newcommand{\refe}[1] {{(\ref {#1})}}%Nov 5
\newcommand{\reft}[1] {{{\textbf{Table} \ref {#1}}}}
\def\la      {{\lambda}}

\def\pmb#1{\setbox0=\hbox{$#1$}%
   \kern-.025em\copy0\kern-\wd0
   \kern.05em\copy0\kern-\wd0
   \kern-.025em\raise.0433em\box0 }

\def\eop{{\hfill $\blacksquare$}}%17Apr01
\def\r{{\rho}}
\newtheorem{theorem}{Theorem}[section]%17Apr01
\newtheorem{lemma}[theorem]{Lemma}%17Apr01
\def\dx     {{\Delta x}}
\def\dz     {{\Delta z}}
\def\dt     {{\Delta t}}

\usepackage[draft]{changes}
\definechangesauthor[name={Wenlong Jin},color=red]{WJ}

\usepackage{multirow}
\usepackage{subcaption}
\usepackage{float}
\usepackage[colorinlistoftodos]{todonotes}

\begin{document}
	\title{Generalized bathtub model of network trip flows} %20190909
\author{Wen-Long Jin \footnote{Department of Civil and Environmental Engineering, California Institute for Telecommunications and Information Technology, Institute of Transportation Studies, 4000 Anteater Instruction and Research Bldg, University of California, Irvine, CA 92697-3600. Tel: 949-824-1672. Fax: 949-824-8385. Email: wjin@uci.edu. Corresponding author}}
\maketitle
\begin{abstract}
\citep{vickrey1991congestion} proposed a bathtub model for the evolution of trip flows served by privately operated vehicles inside a road network based on three premises: (i) treatment of the road network as a single bathtub; (ii) the speed-density relation at the network level, also known as the network fundamental diagram of vehicular traffic, and (iii) the time-independent negative exponential distribution of trip distances. However, the distributions of trip distances are generally time-dependent in the real world, and  Vickrey's model leads to unreasonable results for other types of trip distance distributions. Thus there is a need to develop a bathtub model with more general trip distance distribution patterns.

In this study, we present a unified framework for modeling network trip flows with general distributions of trip distances, including negative exponential, constant, and regularly sorting trip distances studied in the literature. In addition to tracking the number of active trips as in Vickrey's model, this model also tracks the evolution of the distribution of active trips' remaining distances. We derive four equivalent differential formulations from the network fundamental diagram and the conservation law of trips for the number of active trips with remaining distances not smaller than any value. Then we define and discuss the properties of stationary and gridlock states, derive the integral form of the bathtub model with the characteristic method, obtain equivalent formulations by replacing the time coordinate with the cumulative travel distance, and present two numerical methods to solve the bathtub model based on the differential and integral forms respectively. We further study equivalent formulations and solutions for two special types of distributions of trip distances: time-independent negative exponential or deterministic.  In particular, we present six equivalent conditions for Vickrey's bathtub model to be applicable. Finally we demonstrate that the fundamental diagram and the bathtub model can be extended for multi-commodity trip flows with trips served by mobility service vehicles.

\end{abstract}
 
\section{Introduction}
The basic role of a transportation system is to serve trips of persons and goods through privately operated or mobility service vehicles. With privately operated vehicles, passengers and vehicles form integrated units. In this case, travel demand usually depends on the total number of trips and their distribution in both time and space; and supply is restricted to the network size measured by the total number of lane miles as well as road capacities determined by speed limits, vehicle characteristics, and road geometry. With mobility service vehicles provided by transit agencies and transportation network companies, vehicles and passengers are disintegrated. In this case, mobility service vehicles supply seats to satisfy passengers' demand in trips, and the transportation system supplies road capacities for vehicles' demand in the rights-of-way. The imbalance between such demands and supplies leads to congestion on roads, at mobility service stations, and inside vehicles in many metropolitan areas \citep{downs2004traffic,lam1999study,schaller2018new}.

  There can be two types of flows in a road network: (passenger) trip flows and (vehicular) traffic flows.
 For privately operated vehicles, the number of active trips can be converted to the number of running vehicles with given occupancies, and the trip flows are integrated with the traffic flows. However, when trips are served by mobility service vehicles, the number of vehicles and the corresponding traffic dynamics are substantially impacted by mobility service providers' fleet size management strategies, and trip flows and traffic flows are disintegrated. 
Characterizing and modeling congestion dynamics of the two types of flows has been a core task of transportation research.
 Existing traffic flow models generally describe the evolution of vehicular flows with privately operated vehicles and, therefore, can be considered trip flow models.  Such models can be categorized into three types depending on how a road network is treated. In the most dominant approach, a city is divided into zones that serve as origins and destinations of trips, and zones are connected through a discrete network of roads and intersections (uninterrupted or signalized); then travel times and traffic dynamics on each road and intersection can be described by the link performance function \citep[e.g][]{sheffi1984networks}, Webster's delay function \citep[e.g.][]{roess2010traffic}, kinematic wave theory \citep[e.g.][]{jin2012kinematic}, or car-following models \citep[e.g.][]{hidas2005modelling}. This approach employs a detailed representation of a road network and is suitable for design and operational tasks at local levels. However, due to the limitation in resources for coding road networks, usually only major roads are included, and many minor streets and ingress and egress points along arterial roads have to be omitted. In addition, data for time-dependent origin-destination demand patterns are costly to collect and calibrate. Finally, computational costs are high for trip generation, trip distribution, mode choice, departure time choice, and route choice in discrete networks, and the dynamic traffic assignment problem for route choice is still theoretically open and lacks a converging algorithm \citep{peeta2001dta}.
In the second type, a discrete road network is treated as a continuum two-dimensional region \citep{beckmann1952transportation,ho2006continuum}. This type of models capture the approximate spatial distribution of roads and is suitable for describing one-directional flows of vehicles in a region. Such models are more efficient than the first type in network coding and computation, but the two-dimensional origin-destination demand data are still costly to obtain. In addition, these models need to be adjusted for a region with substantially heterogeneous types of roads such as signalized streets and freeways. 
In the third type, a road network is treated as a single unit with roads of the same traffic conditions; in this sense, a road network with vehicular traffic flows is similar to a bathtub with water flow \citep{vickrey1991congestion,vickrey1994types,arnott2013bathtub}. In such bathtub models of network traffic flows, the detailed topology of a network or the spatial distribution of origins, destinations, routes and roads are ignored, and only the lane-miles as well as the average speed-density relation (i.e., network fundamental diagram) are relevant. In addition, time-dependent travel demands are characterized by time-dependent rates of entering trips and the distributions of their distances. Such models are much simpler and directly capture the interactions of dynamic travel demand and network supply even more efficiently. But they have to be modified when origin-destination demands and traffic conditions vary with modes, vehicles, and locations. Note that bathtub models are different from the point queue model with one origin, one destination, and one bottleneck in \citep{vickrey1969congestion}, since there are many origins, destinations, roads, and bottlenecks in a network, even though they are undifferentiated.

Network fundamental diagram, i.e., the functional relationship between the average speed and the average per-lane vehicle density in a road network, was proposed and empirically calibrated for the town center road network in the city of Ipswich, England in \citep{godfrey1969mechanism}. With the availability of more data from loop detectors and other sources, more empirical evidences are obtained recently for other city or freeway networks \citep[e.g.][]{geroliminis2008eus,cassidy2011macroscopic,wada2015empirical,ambuhl2016empirical}. The network fundamental diagram was developed and calibrated for multi-modal traffic and passenger flows in \citep{geroliminis2014three,chiabaut2015evaluation}. However, empirically, it has been shown that links and traffic conditions have to be homogeneous for such a speed-density relation to be well-defined \citep{buisson2009exploring,cassidy2011macroscopic}.  Theoretically, a trapezoidal flow-density relation exists in a signalized network under stationary (periodic in both time and space) traffic conditions  \citep{jin2015performance}. Therefore, such a network-level speed-density relation only captures the static characteristics of traffic flow approximately and should be used with caution when a network contains  heterogeneous types of roads, directional congestion patterns, and highly dynamic traffic conditions. In addition to vehicle densities, the distribution of trips can also impact the average speed of mobility service vehicles due to the boarding and alighting times and other delays; and existing network fundamental diagrams usually ignore the impacts of dynamic trip flows. However, it is a good starting point for studying traffic dynamics at the aggregate level for signalized and freeway networks.

\citet{vickrey1991congestion,vickrey1994types} introduced the first bathtub model for network trip flows served by privately operated vehicles based on three premises: (i) ``a maze of congested streets is treated as an undifferentiated movement area''; (ii) ``movement takes place at a speed which is a function of the density of cars in the area''; and (iii) the trip distance follows a time-independent negative exponential distribution, or the average remaining distance of active trips is constant. As a result, the evolution of the number of active trips is described by an ordinary differential equation. Further, marginal cost pricing was determined based on the instantaneous travel time when a vehicle enters the network. In \citep{small2003hypercongestion}, Vickrey's bathtub model\footnote{In \citep{agnew1976dynamic}, an ordinary differential equation similar to Vickrey's bathtub model was used to model and control a congestion-prone system; but neither of the three premises were mentioned. In \citep{mahmassani1984dynamic}, an ordinary differential equation similar to Vickrey's bathtub model was proposed for traffic flow on a link; but as pointed in \citep{newell1988traffic}, such a model is not applicable for such a single bottleneck. In \citep{daganzo2007gridlock}, the same ordinary differential equation as Vickrey's bathtub model was independently derived based on the first two premises; but the third one was only implicitly assumed. Therefore, we refer to this model as Vickrey's bathtub model.} was used to find the departure time user equilibrium in a network, but it was pointed out that a trip's instantaneous travel time at its entrance is substantially different from its experienced travel time, which depends on the vehicle densities and travel speeds during the whole trip duration.  In \citep{fosgerau2015congestion}, a bathtub model was developed for deterministic distributions of trip distances, and trips are ``regularly sorted'' such that shorter trips enter the network later but exit earlier than longer ones (last-in-first-out); then the departure time user equilibrium was solved with experienced travel times. Yet another bathtub model was developed in \citep{arnott2016equilibrium,arnott2018solving}, in which all trips are assumed to have the same distance, and the dynamics of the number of active trips and the experienced travel times are described by delay-differential equations, which are solved by an iterative method.

This study is motivated by several factors. Empirically, the distribution of trip distances is neither time-independent nor exponential \citep{liu2012understanding,thomas2013empirical}. That is, the assumption of a time-independent exponential distribution of trip distances in Vickrey's bathtub model needs to be relaxed.
Theoretically, many studies following \citep{daganzo2007gridlock} were not aware of the third premise of Vickrey's bathtub model in \citep{vickrey1991congestion} and have attempted to apply Vickrey's bathtub model for constant and other distributions of trip distances; but this has led to physically unreasonable results as information travels too fast   (see \citep{mariotte2017macroscopic} and references therein).
On the other hand, it has been argued that Vickrey's bathtub model is improper as ``the exit rate from downtown traffic (i.e., the arrival rate at work) depends only on the density of downtown traffic, and hence that the first exit occurs as soon as the first entry'' in \citep{arnott2016equilibrium,arnott2018solving}. 
In addition, there seems to be no bathtub model for trip flows served by mobility service vehicles in the literature; but such a model is essential for understanding congestion dynamics in a shared mobility system.
Thus there is a need to develop a bathtub model with more general distance distribution patterns of trips served by both privately operated and mobility service vehicles  and better understand the properties of Vickrey's and other bathtub models in the literature.

In this study we fill the gap by presenting a generalized bathtub model of network trip flows and discussing its properties and solutions, especially its relationship with existing bathtub models in the literature. First we assume that a trip is served by a privately operated vehicles,  and we use trips and vehicles interchangeably in this case.
This model generalizes Vickrey's by relaxing the third premise of the latter. That is, this model works for any distributions of trip distances. In Vickrey's bathtub model, the main state variable is the total number of active trips; in contrast, in the new model, the main state variable is the number of active trips with remaining distances not smaller than a value. From the conservation of trips, we derive four equivalent partial differential equations to track the evolution of variables related to the distribution of remaining distances. Then we discuss the properties as well as analytical and numerical solutions of the model. In particular, we show that Vickrey's bathtub model is a special case when the trip distance follows a time-independent negative exponential distribution.  Finally we extend the bathtub model for trips served by mobility service vehicles.
For privately operated vehicles,  the bathtub model also describes the dynamics of both trip and traffic flows. However, when trips are served by mobility service vehicles, the bathtub model only describes the dynamic progression of trips in a road network, and the vehicular traffic flow dynamics have to be developed separately. Therefore, the bathtub model is for network trip flows, not necessarily for network traffic flows.

The rest of the article is organized as follows. In Section 2, after defining the number of active trips at a time instant, that with a remaining distance not smaller than a value, the cumulative travel distance, and other variables, we derive four equivalent formulations of the model from the conservation of trips. In Section 3, we define and discuss the properties of stationary and gridlock states, derive the integral form of the bathtub model with the characteristic method, obtain equivalent formulations in different coordinates, and present a numerical method and an example. In Sections 4 and 5, we further studied equivalent formulations and solutions for two special types of distributions of trip distances: time-independent negative exponential or deterministic. In Section 6, we extend the model for trips served by mobility service vehicles and multi-commodity trip flows. In Section 7, we conclude the study with discussions.

\section{Definitions of variables and derivation of the model}

A list of notations is given in \reft{table:notations}.

\btb\bc
\begin{tabular}{|c|l||c|l|}\hline
	Variables & Definitions& Variables &Definitions \\\hline
	$t$ & Time &  $x$ & Trip distance \\\hline
	$\dt$& Time interval &$\dx$& Distance interval \\\hline\hline
	$L$ & Total lane-miles of a network & $v(t)$ & Average travel speed  \\\hline
	$z(t)$ & Cumulative travel distance & $\tau(z)$& Cumulative travel time \\\hline \hline
	$\tilde \varphi(t,x)$ & \multicolumn{3}{l|} {Probability density function of the entering trip distance at $t$} \\\hline
	$\tilde \Phi(t,x)$ & \multicolumn{3}{l|} {Proportion of the entering trips with distances not smaller than $x$}\\\hline 
		$\tilde B(t)$ & \multicolumn{3}{l|} {Average entering trips' distances at $t$} \\\hline
	$\Upsilon(t,x)$ & \multicolumn{3}{l|} {Travel time for a trip entering at $t$ with a distance of $x$ }\\\hline
	$\bar \Upsilon(t)$ & \multicolumn{3}{l|} {Average travel time for all trips entering at $t$}\\\hline	
	$\theta(t,x)$ & \multicolumn{3}{l|} {Effective distance of a trip entering at $t$ with a distance of $x$ }\\\hline \hline
	$F(t)$ & In-flow (entering trip flow)& $G(t)$ & Out-flow (exiting trip flow) \\\hline
	$f(t)$ & In-flux (entering trip rate) & $g(t)$ & Out-flux (exiting trip rate) \\\hline
	$\lambda(t)$& Number of active trips at $t$ & $\r(t)$   & Per-lane density \\\hline
	$\varphi(t,x)$ & \multicolumn{3}{l|} {Probability density function of the active trips' remaining distances at $t$}  \\\hline
	$\Phi(t,x)$& \multicolumn{3}{l|} {Proportion of the trips with remaining distances not smaller than $x$}\\\hline
	$k(t,x)$ & \multicolumn{3}{l|} {Density of active trips at $t$ with a remaining distance $x$}  \\\hline 
	$K(t,x)$ & \multicolumn{3}{l|} {Number of active trips at $t$ with a remaining distance not smaller than $x$ }\\\hline
	 $B(t)$  & \multicolumn{3}{l|} {Average remaining distance of active trips}  \\\hline\hline
	 	$n$ & Trip $n$  & $X(t,n)$ & Remaining distance of trip $n$ at $t$  \\\hline
	 	$T(n,x)$ & \multicolumn{3}{l|} {Time for trip $n$'s remaining distance is $x$} \\\hline 
	 	$N(t,x)$ & \multicolumn{3}{l|} {Cumulative flows of trips at $(t,x)$ } \\\hline

\end{tabular}
\ec
\caption{List of notations} \label{table:notations}
\etb

\subsection{Definitions of variables}

\bfg\bc
\includegraphics[width=5in]{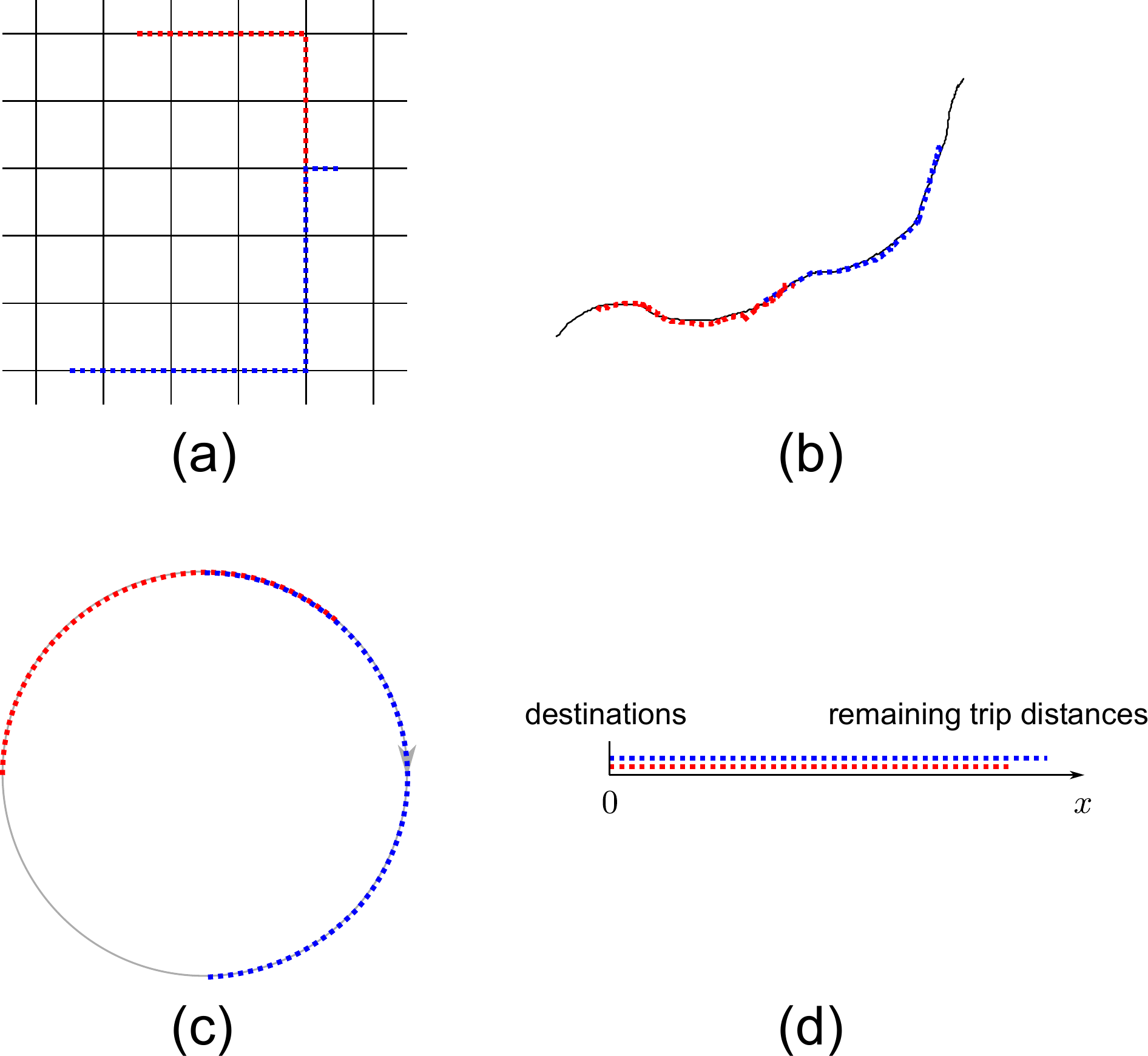}
\caption{Illustration of traffic networks with many origins, destinations, and bottlenecks: (a) A signalized city road network; (b) A freeway/arterial corridor; (c) A beltway or a two-way bus/metro line; (d) Trips in the $x$ space. }\label{city_traffic_network}
\ec\efg

For the three examples of road networks in \reff{city_traffic_network}(a)-(c), their sizes are measured by the lane-miles, $L$.
We assume that  there are many origins and destinations, the capacity of one location is shared by many traffic streams, and there can be many bottlenecks in the network. In particular, the beltway network in  \reff{city_traffic_network}(c) can be considered a limit case of the beltway considered in \citep{daganzo1996gridlock} but with many on- and off-ramps, or two-way bus or metro lines with many stations. 

We denote the number of active trips on the network at $t$ by $\lambda(t)$ (unit: number of trips), which equals the number of running vehicles when each trip is served by one privately operated vehicle. Then the vehicle density per lane equals $\r(t)$ (unit: vehicles per lane-mile):
\bqn
\r(t)=\frac{\la(t)}L. \label{la2rho}
\eqn

The average travel speed is denoted by $v(t)$ (unit: mph). 
Here we assume that all roads in a network have the same average traffic conditions at any moment, and we can treat the whole network as a single unit, or a bathtub.\footnote{Note that, even when traffic conditions are periodic both in time and space in a signalized network, the traffic conditions, including both densities and speeds, vary with time and location. Thus, the homogeneity in traffic conditions should be understood at the average sense over a relatively long period of time or at the statistical sense.} 
We denote the cumulative travel distance by $z(t)$; i.e., $z(t)=\int_0^t v(s) ds$, or equivalently $z(0)=0$, and 
\bqn
\dot z(t)= v(t). \label{def:zt}
\eqn
When vehicles always move with $v(t)>0$, we denote the inverse function of $z(t)$ as $\tau(z)$,
\bqn
\tau(z(t))&=&t,
\eqn
 which is the cumulative travel time. Thus
\bqn
\dot \tau(z)&=& \frac 1 {v(z)}, \label{def:tau}
\eqn
where $v(z)=v(\tau(z))$.

During a study period during $[0,T]$, each trip has an origin, a destination, a path, an entering time,  and an exiting time. In the bathtub model, a trip can be represented by its entering time, $t$, and its distance, $x$. We denote the effective distance of a trip entering at $t$ with a distance of $x$ by $\theta(t,x)=x+z(t)$, which is the cumulative travel distance when the trip exits the network. In addition, the exit time for this trip is effectively the same as a trip that is in the network at $t=0$ with a distance of $\theta(t,x)$. Since $z(t)$ increases in $t$ with positive speeds, trips with shorter effective distances exit the network earlier. We refer to this as the shorter-(effective) distance-first-out principle. 

We define the rate of entering trips in the whole network  by $f(t)$ (unit: vph), which is also the in-flux. \footnote{For privately operated vehicles, the in-flux is related to the traffic queueing process in a parking lot, side streets not included in the road network, or roads outside the network. It is related to the queues of passengers waiting to board buses or shared cars in a shared mobility system.} Then $\frac 1 L f(t)$ is the number of entering trips per unit time and unit lane-mile. 
\footnote{As an example, consider a grid network comprised of $\beta$ blocks, in which each link's length is $b$, and the number of lanes is $l$.  Then the total lane-miles is $L \approx 2 \beta  b l$, since each link  is shared by two blocks, and the area is $A=\beta b^2$. Assume that the number of entering trips is proportional to the area, and the trip generation rate at $t$ is $a(t)$ (unit: trips/mile$^2$/hr). In this case, the in-flux per lane-mile is $\frac 1L f(t)= \frac{A a(t)}{L}=\frac{1}{2 l} b a(t)$,
which is independent of the number of blocks, but increases with the link length, and decreases in the road width. This simple model can be used to study the impacts of housing and job distributions on the in-flux.}
We denote the cumulative in-flow by $F(t)$: $F(t)=\int_0^t f(s)ds$.

We denote the probability density function of entering trips' total distances, $x$, at $t$ by $\tilde \varphi(t,x)$ (unit: mile$^{-1}$), which satisfies the following conditions: (i) $\tilde \varphi(t,x) \geq 0$; (ii) $\int_0^ \infty \tilde \varphi(t,x) dx =1$; (iii) $\tilde \varphi(t,\infty)=0$. 
We denote the proportion of the entering trips with distances not smaller than $x$ by $\tilde \Phi(t,x)$: $\tilde \Phi(t,x)=\int_x^\infty \tilde \varphi(t,y)dy$, 
and $\tilde \varphi(t,x)=-\pd{}x \tilde \Phi(t,x)$.
Thus, $\tilde \Phi(t,0)=1$, $\tilde \Phi(t,\infty)=0$, and $\tilde \Phi(t,x)$ is non-increasing in $x$.
If we denote the average distance of the entering trips at time $t$ by $\tilde B(t)$, then 
\bqn
\tilde B(t)&=&\int_0^\infty x \tilde \varphi(t,x) dx=\int_0^\infty \tilde \Phi(t,x) dx. \label{averagetriplength}
\eqn
Here $f(t)$ and $\tilde \Phi(t,x)$ represent the demand pattern of trips.

For a trip with an entering time of $t$ and  a distance of $x$, we denote its travel time by $\Upsilon(t,x)$:
\bqn
\Upsilon(t,x)&=&\tau(x+z(t))-t. \label{def:triptraveltime}
\eqn
We denote the average travel time for entering trips at $t$ by $\bar \Upsilon(t)$. Then we have
\bqn
\bar \Upsilon(t)&=&\int_0^\infty \Upsilon(t,x) \tilde \varphi(t,x) dx. \label{def:traveltime1}
\eqn
It can be easily shown that
\bqn
\bar \Upsilon(t)&=&\int_0^\infty \tilde \Phi(t,x) \frac 1{ v(\tau (x+z(t)) )} dx, \label{def:traveltime}
\eqn
which can be approximated by
\bsq
\bqn
\bar \Upsilon(t)&\approx& \frac{\tilde B(t)}{v(t)}, \label{def:traveltime:approx1}
\eqn
as in \citep{vickrey1991congestion},
or
\bqn
\bar \Upsilon(t)&\approx& \frac{\tilde B(t)}{v(\tau(z(t)+\tilde B(t)))}. \label{def:traveltime:approx2}
\eqn
\esq
In \refe{def:traveltime:approx1}, the travel speeds are assumed to be the same as the instantaneous speed when the trips enter the network. In \refe{def:traveltime:approx2}, the travel speeds are assumed to be the same as the instantaneous speed when the trip with the average distance exits the network.

A network can be considered as a trip processing machine: trips with a time-dependent trip distance distribution of $\tilde \varphi(t,x)$ are added into it dynamically at a rate of $f(t)$, and each trip's remaining distance is shredded by $v(t)$ per unit time until it exits the network. At time $t$, an active trip in the network has a remaining distance $x>0$, and different trips can be represented by a (dotted) line segment in the $x$ space, as shown in \reff{city_traffic_network}(d). 
We denote the probability density function of the remaining trip distance, $x$, at $t$ by $\varphi(t,x)$ (unit: mile$^{-1}$), which satisfies the following conditions: (i) $\varphi(t,x) \geq 0$; (ii) $\int_0^ \infty \varphi(t,x) dx =1$; (iii) $\varphi(t,\infty)=0$. 
Furthermore, if we denote the proportion of the trips with remaining distances not smaller than $x$ by $\Phi(t,x)$, then $\Phi(t,x)=\int_x^\infty \varphi(t,y)dy$, and $\varphi(t,x)=-\pd{}x \Phi(t,x)$.
Therefore, $\Phi(t,0)=1$, $\Phi(t,\infty)=0$, and $\Phi(t,x)$ is non-increasing in $x$. 
We denote the average remaining trip distance at $t$ by $B(t)$:
\bqn
B(t)&=&\int_0^\infty x \varphi(t,x)dx=\int_0^\infty \Phi(t,x) dx. \label{averageremainingtripdistance}
\eqn

The density of active trips with a remaining distance of $x$ is denoted by $k(t,x)$ (unit: vehicles per mile): 
\bsq
\bqn
k(t,x)&=&\lambda(t) \varphi(t,x),
\eqn
which leads to
\bqn
\lambda(t)&=&\int_0^\infty k(t,x)dx. \label{k2la}
\eqn
\esq
We define the number of trips at $t$ with a remaining distance not smaller than $x$ by $K(t,x)$ (unit: number of vehicles). In this study this will be the main variable.  Note that the variable of $K(t,x)$ was first defined in \citep{vickrey1991congestion}; but its evolution dynamics were not explicitly considered. Then
\bsq
\bqn
K(t,x)&=& \la (t) \Phi(t,x).
\eqn
That is, if $K(t,x)=K$, then the $K$th longest active trip has a remaining distance of $x$ at $t$. 
Thus $\pd {}x K(t,x)=- k(t,x)$, $K(t,x)=\int_x^\infty k(t,y) dy$, and  
\bqn
\la (t)&=&K(t,0), \label{K2la}
\eqn
\esq
since trips with negative distances are assumed to have exited the network.

From $t$ to $t+\dt$ for a small time interval $\dt$, the number of exiting (completing) trips is $g(t)\dt$, which equals those with a remaining distance between $0$ and $v(t) \dt$ at time $t$. \footnote{Here we assume that all trips can exit the network once their remaining distances are non-positive. Note that this may not be the case when vehicles have to cruise for parking \citep{shoup2006cruising}.}
That is, 
\bqs
g(t)\dt &=& \int_{0}^{v(t)\dt} k(t,x) dx\approx k(t,0) v(t) \dt; \label{discrete-g}
\eqs
or equivalently, the out-flux (the completion rate of trips) equals
\bqn
g(t) &=& k(t,0) v(t)=\varphi(t,0) \lambda(t) v(t) , \label{continuous-g}
\eqn
which depends on the number of active trips, the travel speed, and the probability density function of the remaining distance $0$.
This can be considered a generalized version of the network exit function introduced in \citep{gonzales2012morning}. The cumulative out-flow is denoted by $G(t)=\int_0^t g(s) ds$.

\subsection{Network fundamental diagram}

By the definition of the network fundamental diagram \citep{godfrey1969mechanism,vickrey1991congestion}, we assume the following network-level speed-density relation when a trip is served by a privately operated vehicle: 
\bqn
v(t)&=&V(\r(t))=V\left(\frac{\la(t)}L\right), \label{speed-density-relation}
\eqn
which represents the supply of a road network.
For freeway networks, we can use the triangular fundamental diagram \citep{cassidy2011macroscopic}:
\bqn
V(\r)=\min\{u, w (\frac {\kappa} \r-1)\}, \label{triangular-fd}
\eqn
where $u$ is the free-flow speed, $w$ the shock wave speed in congested traffic, and $\kappa$ the per-lane jam density.
For signalized road networks, we can use the following trapezoidal fundamental diagram \citep{jin2015performance}:
\bqn
V(\r)&=&\min\{u, \frac C{\rho}, w (\frac {\kappa} \r-1) \},
\eqn
where the additional parameter $C$ is the average capacity. 
In the literature, the Greenshields fundamental diagram is also used:
\bqn
V(\r)&=& u (1-\frac{\r}\kappa).
\eqn
For freeway networks, $u$ and $w$ depend on speed limits and vehicle characteristics; for signalized networks, $u$, $w$, and $C$ also depend on signal settings.
From the speed-density relation, a trip entering the network at a later time could impact an earlier one still in the network as the former changes the number of active trips.

The corresponding average per-lane flow-density relation can be written as
\bqn
q&=&Q(\r)\equiv \r V(\r).
\eqn
Generally, the travel speed is non-increasing in the per-lane density; i.e., $\der {}\r V(\r)\leq 0$; and the flow-density relation is concave. \footnote{Sometimes the travel speed can increase in the density as in night traffic \citep{leveque2001night}. Non-concave flow-density relations have also be reported in the literature \citep[e.g.][]{kerner1993cluster}.}  The model developed in this study works for very general speed-density and flow-density relations. For examples, the speed-density relation can be piecewise constant, and the flow-density relation can be non-concave or discontinuous.

When $\r(t)=\kappa$ and $\la(t)=L\kappa$, all roads in the network are jammed with vehicles, and both the travel speed and the flow-rate are zero. In this case, the network reaches a standstill or gridlock state \citep{vickrey1969congestion,daganzo2007gridlock}.

\subsection{Three conservation laws and four equivalent differential formulations of the generalized bathtub model}
At $t$, the trip processing rate is given by
\bqn
\la(t) v(t)&=& L Q\left(\frac {\lambda(t)}L\right).
\eqn 
Therefore the maximum trip processing rate equals $LC$, where $C=\max_\r Q(\r)$. 
For a concave flow-density relation, the trip processing rate decreases with increasing vehicle densities in hypercongestion \citep{small2003hypercongestion}. 
From the conservation of trip-miles we have
\bqn
\la(0) B(0)+\int_0^t f(s) \tilde B(s) ds -\int_0^t \la(s) v(s) ds &=&\la(t) B(t), \label{tri-miles-conservation}
\eqn
where $\la(0) B(0)$ is the initial trip-miles, $\int_0^t f(s) \tilde B(s) ds$ the added trip-miles till $t$, $\int_0^t \la(s) v(s) ds$ processed trip-miles till $t$, and $\la(t) B(t)$ the remaining trip-miles. 

In addition, from the conservation of total trips, we can calculate the cumulative out-flow as
\bqn
G(t)&=&\la(0)+F(t)-\la(t), \label{conservation-total-trips}
\eqn
whose differential version is $g(t)=f(t)-\dot \la(t)$, or
\bqn
\dot \la(t)&=&f(t)-g(t). \label{diff-conservation-total-trips}
\eqn
From \refe{continuous-g} and \refe{diff-conservation-total-trips}, we have the following evolution equation for the total number of active trips
\bqn
\dot \la (t)&=& f(t) - \varphi(t,0) \la(t) v(t). \label{continuous-la-eqn}
\eqn
Incorporating the network speed-density relation in \refe{speed-density-relation}, the evolution of the number of active trips in the network is governed by the following ordinary differential equation:
\bqn
\dot \la (t)&=& f(t) - \varphi(t,0) \la(t) V\left(\frac{\la(t)}L\right). \label{la-ode}
\eqn
In this equation, $\la(t)$ is the unknown variable, $f(t)$ and $V(\cdot)$ are given demand and supply. Note that \refe{la-ode} is an incomplete model, as $\varphi(t,0)$ is also unknown. 

We further complete \refe{la-ode} by tracking the probability density function of the remaining trip distance, including $\varphi(t,0)$. At $t+\dt$, the number of active trips with a remaining distance between $x$ and $x+\dx$ is $\int_{x}^{x+\dx} \la (t+\dt) \varphi(t+\dt, y) dy= \la(t+\dt) \varphi(t+\dt,x) \dx$. Such trips include (i) those with a remaining distance between $x+v(t)\dt$ and $x+v(t)\dt+\dx$ at $t$, which will be shortened to be between $x$ and $x+\dx$ at $t+\dt$, and (ii) those entering the network between $t$ and $t+\dt$ with a trip distance between $x$ and $x+\dx$. The number of the first group of trips is $\la(t) \varphi(t,x+v(t)\dt) \dx$; and the number of the second group is $f(t) \dt \tilde \varphi(t,x) \dx$. Therefore we have (with $\dx$ eliminated from all terms)
\bqn
\la (t+\dt) \varphi(t+\dt,x)&=& \la (t) \varphi(t,x+v(t)\dt) + f(t) \dt \tilde \varphi(t,x). \label{discrete-pdf-eqn}
\eqn
Replacing $\la(t+\dt)$ by $\la(t)+\dot\la(t) \dt$ and $\varphi(t,x+v(t)\dt)$ by $\varphi(t,x)+\varphi_x(t,x) v(t)\dt$, we obtain from the above equation
\bqs
\la(t)[\varphi(t+\dt,x) -\varphi(t,x)-v(t) \varphi_x(t,x) \dt ]&=& f(t) \tilde \varphi(t,x) \dt-\dot \la(t) \varphi(t+\dt,x) \dt.
\eqs
Dividing both sides of the above equation by $\dt$ and let $\dt\to 0$, we obtain the following continuous version of \refe{discrete-pdf-eqn}:
\bqs
\varphi_t(t,x) -v(t) \varphi_x (t,x)&=&\frac 1 {\la(t)} [ f(t) \tilde \varphi(t,x) -  \varphi(t,x) \dot \la(t) ].
\eqs
Further by substituting $\dot \la(t)$ in \refe{continuous-la-eqn} into the above equation, we have 
\bqn
\pd{}t \varphi(t,x) -v(t) \pd{}x \varphi(t,x)&=&\frac{f(t)}{\la(t)} (\tilde \varphi(t,x) -\varphi(t,x))+ \varphi(t,x)\varphi(t,0) v(t). \label{continuous-pdf-eqn}
\eqn
Thus  \refe{speed-density-relation}, \refe{continuous-la-eqn}, and \refe{continuous-pdf-eqn} form a complete bathtub model, in which $\la(t)$, $v(t)$, and $\varphi(t,x)$ are the unknown variables. Hereafter we refer to this model as the $\varphi$-model.

Integrating both sides of \refe{continuous-pdf-eqn} with respect to $x$ from $x$ to $\infty$, we obtain the following equation of $\Phi(t,x)$:
\bqn
\pd{}t \Phi(t,x) -v(t) \pd{}x \Phi(t,x)&=&\frac{f(t)}{\la(t)} (\tilde \Phi(t,x) -\Phi(t,x))+ \Phi(t,x)\varphi(t,0) v(t), \label{continuous-cdf-eqn}
\eqn
which, with \refe{speed-density-relation}, \refe{continuous-la-eqn}, and $\varphi(t,0)=-\Phi_x(t,x) \big|_{x=0}$, forms the $\Phi$-model. In this model, $\la(t)$, $v(t)$, and $\Phi(t,x)$ are the unknown variables.
The corresponding discrete $\Phi$-model is
\bqn
\la (t+\dt) \Phi(t+\dt,x)&=& \la (t) \Phi(t,x+v(t)\dt) + f(t) \dt \tilde \Phi(t,x). \label{discrete-cdf-eqn}
\eqn

An equivalent discrete evolution equation of \refe{discrete-pdf-eqn} is 
\bqn
k(t+\dt,x)&=&k(t,x+v(t)\dt)+f(t) \tilde \varphi(t,x) \dt. \label{discrete-k-eqn}
\eqn
Replacing $k(t,x+v(t)\dt)$ by $k(t,x)+\pd{}x k(t,x) v(t)\dt$ in the above equation, we obtain its continuous version as
\bqn
\pd {}t k(t,x) -v(t) \pd {}x k(t,x)&=& f(t) \tilde \varphi(t,x). \label{continuous-k-eqn}
\eqn
Thus \refe{k2la}, \refe{speed-density-relation}, and \refe{continuous-k-eqn} yield another formulation of the bathtub model, in which $\la(t)$, $v(t)$, and $k(t,x)$ are the unknown variables. Plugging \refe{k2la} and \refe{speed-density-relation} into \refe{continuous-k-eqn}, we obtain the following simplified integral-differential equation:
\bqn
\pd {}t k(t,x) -V\left(\frac 1L \int_0^\infty k(t,y) dy \right) \pd {}x k(t,x)&=& f(t) \tilde \varphi(t,x), \label{k-model}
\eqn
where $k(t,x)$ is the only unknown variable. Hereafter we refer to \refe{k-model} or, equivalently, \refe{k2la}, \refe{speed-density-relation}, and \refe{continuous-k-eqn}, as the $k$-model.

\bfg\bc
\includegraphics[width=5in]{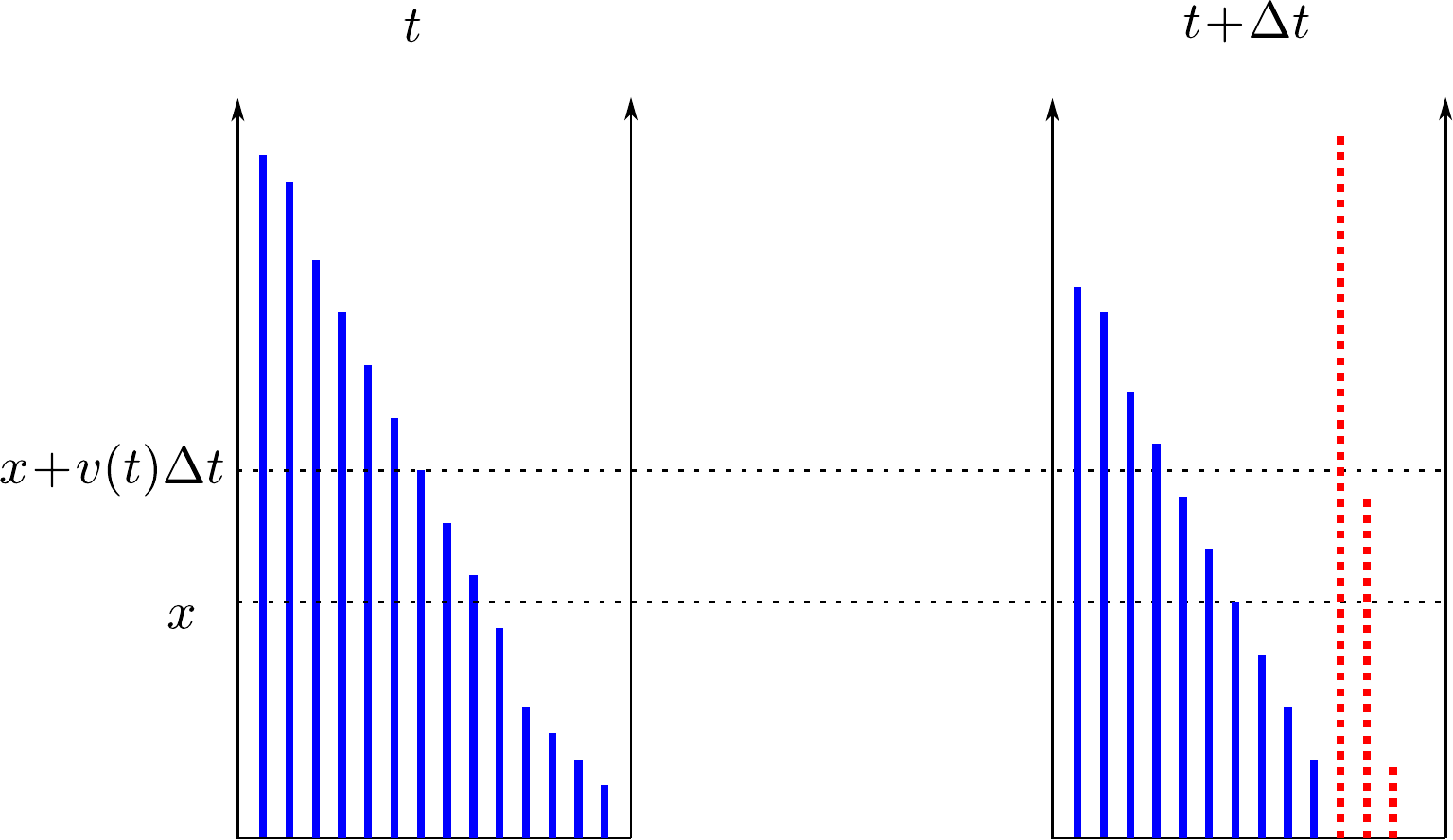}
\caption{Illustration of the discrete $K$-model}\label{rho-model-illustration}
\ec\efg

Integrating both sides of \refe{continuous-k-eqn} from $x$ to $\infty$ with respect to $x$, we obtain 
\bqn
\pd {}t K(t,x) -v(t) \pd {}x K(t,x)&=& f(t) \tilde \Phi(t,x), \label{continuous-K-eqn}
\eqn
which can also be derived from the following discrete version:
\bqn
K(t+\dt,x)&=&K(t,x+v(t)\dt)+f(t) \tilde \Phi(t,x) \dt. \label{discrete-K-eqn}
\eqn
An illustration of the derivation process of the above equation is in \reff{rho-model-illustration}, in which the left and right sides represent the traffic conditions at $t$ and $t+\dt$ respectively, the solid blue lines represent the trips in the network at $t$,  the  dotted red lines represent the trips entering the network between $t$ and $t+\dt$. From the figure, we can see that the number of trips with a remaining distance not smaller than than $x$ at $t+\dt$ (the left seven solid blue lines and the left two dotted red lines in the right figure), $K(t+\dt,x)$ in \refe{discrete-K-eqn}, equals the number of trips with a remaining distance not smaller than than $x+v(t)\dt$ (the left seven solid blue lines in the left figure), $K(t,x+v(t)\dt)$  in \refe{discrete-K-eqn}, plus the number of entering trips with a distance not smaller than $x$ (the left two dotted red lines in the right figure), $f(t) \tilde \Phi(t,x) \dt$ in \refe{discrete-K-eqn}.
From the definitions of $\Phi(t,x)$ and $K(t,x)$, we can see that \refe{continuous-la-eqn} is a special case of \refe{continuous-K-eqn} when $x=0$.
Thus \refe{K2la}, \refe{speed-density-relation}, and \refe{continuous-K-eqn} yield another formulation of the model, in which $\la(t)$, $v(t)$, and $K(t,x)$ are the unknown variables. Plugging \refe{K2la} and \refe{speed-density-relation} into \refe{continuous-K-eqn}, we obtain the following simplified partial differential equation:
\bqn
\pd {}t K(t,x) -V\left(\frac {K(t,0)}L \right) \pd {}x K(t,x)&=& f(t) \tilde \Phi(t,x), \label{K-model}
\eqn
where $K(t,x)$ is the only unknown variable. Hereafter we refer to \refe{K-model} or, equivalently, \refe{K2la}, \refe{speed-density-relation}, and \refe{continuous-K-eqn}, as the $K$-model.

We can see that the generalized bathtub model is derived from the network fundamental diagram and the conservation law for the number of active trips with remaining distances not smaller than any value. Note that the conservation law for the total number of active trips is a special case of this one.

\section{Properties and solutions}

In this section we discuss the properties and solutions of the generalized bathtub model, in particular, of the $K$-model.

\subsection{Stationary states and their stability} \label{section:stationarystate}

We define a traffic state as stationary if for any $x$ the state variable is time-independent;
\bqn
\pd{}t K(t,x)&=&0. \label{def:stationarystate}
\eqn

\begin{lemma} The traffic state is stationary if and only if
	the number of active trips is constant: $\la(t)=\la$; the travel speed is constant: $v(t)=V(\frac {\la}L)=v$;  the remaining trip distance distribution is time-independent: $\Phi(t,x)=\Phi(x)$; and the average remaining trip distance $B(t)=B$ is constant.
\end{lemma}
{\em Proof}. 
Setting $x=0$ in \refe{def:stationarystate}, we have $\dot \la(t)=0$; 
i.e., the number of active trips is constant: $\la(t)=\la$. Hence the travel speed $v(t)=V(\frac {\la}L)=v$ is constant. Further from \refe{def:stationarystate} we have $\pd{}t \Phi(t,x)=0$. 
Thus $\Phi(t,x)=\Phi(x)$ is time-independent. From the definition of $B(t)$, it is also constant.  Hence the necessary conditions are proved. The sufficient conditions are obvious from the definition of $K(t,x)$. \eop

\begin{theorem} \label{thm:stasta}

In a stationary state of the generalized bathtub model, the in-flux is constant
\bqn
f(t)&=&f\equiv L Q\left(\frac {\la}L\right) \varphi(0),
\eqn
and the trip distance distribution is time-independent
\bqn
\tilde \Phi(t,x)&=&\tilde \Phi(x)\equiv\frac{\varphi(x)}{\varphi(0)}.
\eqn
The average trip distance is constant at
\bqn
\tilde B(t)&=&\tilde B\equiv \frac 1{\varphi(0)},
\eqn
which may not equal $B$.
\end{theorem}

{\em Proof}. In a stationary state, \refe{K-model} can be simplified as
\bqs
 L Q\left(\frac {\la}L\right) \varphi(x)&=&f(t) \tilde \Phi(t,x),
\eqs
which at $x=0$ leads to a constant in-flux: $f(t)=f\equiv L Q\left(\frac {\la}L\right) \varphi(0)$,
and a time-independent trip distance distribution $\tilde \Phi(t,x)=\frac{\varphi(x)}{\varphi(0)}$.
Then the average trip distance is constant $\tilde B(t)=\tilde B\equiv \frac 1{\varphi(0)}$.
\eop

Near a stationary state with $K(t,x)=\la_0 \Phi(x)$ and $f(t)\tilde \Phi(t,x)=LQ\left(\frac {\la_0}L\right) \varphi(x)$, if a small disturbance $\epsilon(t)$ is applied to the initial state, such that \refe{la-ode} can be written as
\bqn
\dot \epsilon(t)&=&LQ\left(\frac {\la_0}L\right) \varphi(0)-LQ\left(\frac {\la_0+\epsilon(t)}L\right) \varphi(0) \approx - L \dot Q\left(\frac {\la_0}L\right) \varphi(0) \epsilon(t),
\eqn
which is stable if $\dot Q\left(\frac {\la_0}L\right) \geq 0$ and unstable if $\dot Q\left(\frac {\la_0}L\right) < 0$. An unstable system settles down at a local minimum state of $Q(\r)$. In particular, for a concave flow-density relation, a stationary state in hypercongestion is unstable, and the system eventually converges to a gridlock state with a small disturbance.

If a network reaches a gridlock state at $t$; i.e., $\la(t)=L\kappa$, then the out-flux $g(t)=0$, and \refe{la-ode} leads to $\dot \la(t)=f(t)\geq 0$ and $\la(t) \geq L\kappa$ after $t$. Such a gridlock state usually occurs under non-recurrent, abnormal conditions, as people can choose their departure times to avoid gridlock under normal conditions.
\refe{tri-miles-conservation} can be used to determine whether gridlock develops in a network. For example, if both $f(t)=f$ and $\tilde B(t)=\tilde B$ are time-independent, the network becomes gridlocked  if 
\bqn
f\tilde B> LC. \label{gridlock-sufficientcondition}
\eqn
Here $f\tilde B$ represents the demand, which equals the product of the in-flux and the average trip distance, and $LC$ the supply, which equals the product of the road capacity per-lane and the total lane-miles of the network. 
Thus a network becomes gridlocked when the demand is higher than the supply for a long period of time. 
Note that \refe{gridlock-sufficientcondition} is the sufficient condition for gridlock but may not be necessary, since the initial condition also plays a vital role in the development of gridlock. 
A network cannot be brought out of gridlock by reducing the demand alone. In this case, the only remedy is to divert vehicles away from their intended destinations; i.e., some trips with positive remaining distances are forced to exit the network. For example, if all trips with a remaining distance smaller than $x_0>0$ have to exit, then the out-flux in \refe{discrete-g} is revised as
\bqn
g(t)\dt&=&\int_0^{x_0} k(t,x) dx =\la(t) (1-\Phi(t,x_0)),
\eqn 
which is positive even in the gridlock state.

\subsection{Characteristic curve and the integral bathtub model}\label{section:general:analytical}
We define a characteristic curve in the $(t,x)$ space as $x=x(t)$ with $\dot x(t)=-v(t)$. 
Then a characteristic curve emanated from $(0,x_0)$ can be written as $x(t)=x_0-z(t)$. For a trip entering the network at $t=0$ with a trip distance $x_0$, its trajectory is also given by $x(t)=x_0-z(t)$. Therefore, the characteristic curves coincide with the trip trajectories.

From \refe{K-model} we have an ordinary differential equation for $K(t,x_0-z(t))$ along the characteristic curve:
\bqs
\der{}t K(t,x_0-z(t))&=&\pd{}t  K(t,x_0-z(t)) -v(t) \pd{}x  K(t,x_0-z(t)) =f(t) \tilde \Phi(t,x_0-z(t)), \eqs
from which we can solve $K(t,x_0-z(t))$ as follows:
\bqs
K(t,x_0-z(t))&=&K(0,x_0)+\int_0^t f(s) \tilde \Phi(s, x_0-z(s)) ds,
\eqs
where $K(0,x_0)$ is given by the initial condition, and $\int_0^t f(s) \tilde \Phi(s, x_0-z(s)) ds$ by the boundary condition.
Thus, if $z(t)$ is known, $K(t,x)$ is solved by
\bqn
K(t,x)&=&K(0,x+z(t))+\int_0^t f(s) \tilde \Phi(s, x+z(t)-z(s)) ds, \label{integral-bathtub}
\eqn
which is the integral form of the $K$-model and referred to as the integral $K$-model. \reff{integral-K-conservation}(a) explains \refe{integral-bathtub} graphically: since $K(t,x)$ represents the number of trips with a distance not shorter than $x$ at $t$, it includes the initial trips that are not shorter than $x+z(t)$ and all trips that enter the network at $s\leq t$ with a distance not shorter than $x+z(t)-z(s)$; the number of such initial trips equals $K(0,x+z(t))$, and the number of such entering trips equals  $\int_0^t f(s) \tilde \Phi(s, x+z(t)-z(s)) ds$; and all such trips are in the shaded region in the figure.
From \refe{integral-bathtub}, we have the following integral form for the total number of active trips
\bsq\label{integral-lambda}
\bqn
\la(t)&=&K(t,0)=K(0,z(t))+\int_0^t f(s) \tilde \Phi(s, z(t)-z(s)) ds. 
\eqn
Similarly \reff{integral-K-conservation}(b) explains \refe{integral-lambda}(a) graphically.
\refe{integral-lambda}(a) and \refe{def:zt}, 
\bqn
\dot z(t)=V\left (\frac{\la(t)}L \right ),
\eqn
\esq
 form a system of differential-integral equations, from which $\la(t)$ and $z(t)$ can be solved. Furthermore, $K(t,x)$ can be solved from \refe{integral-bathtub}, and $G(t)$ can be solved from  \refe{conservation-total-trips}. As an example, if the entering trips' distances follow a time-dependent exponential distribution, $\tilde \Phi(t,x)=e^{-\frac{x}{\tilde B(t)}}$, \refe{integral-lambda} can be simplified as
\bqs
\la(t)&=&K(0,z(t))+\int_0^t f(s) e^{-\frac{z(t)-z(s)}{\tilde B(t)}} ds. 
\eqs	
 
 \bfg\bc
 \includegraphics[width=5in]{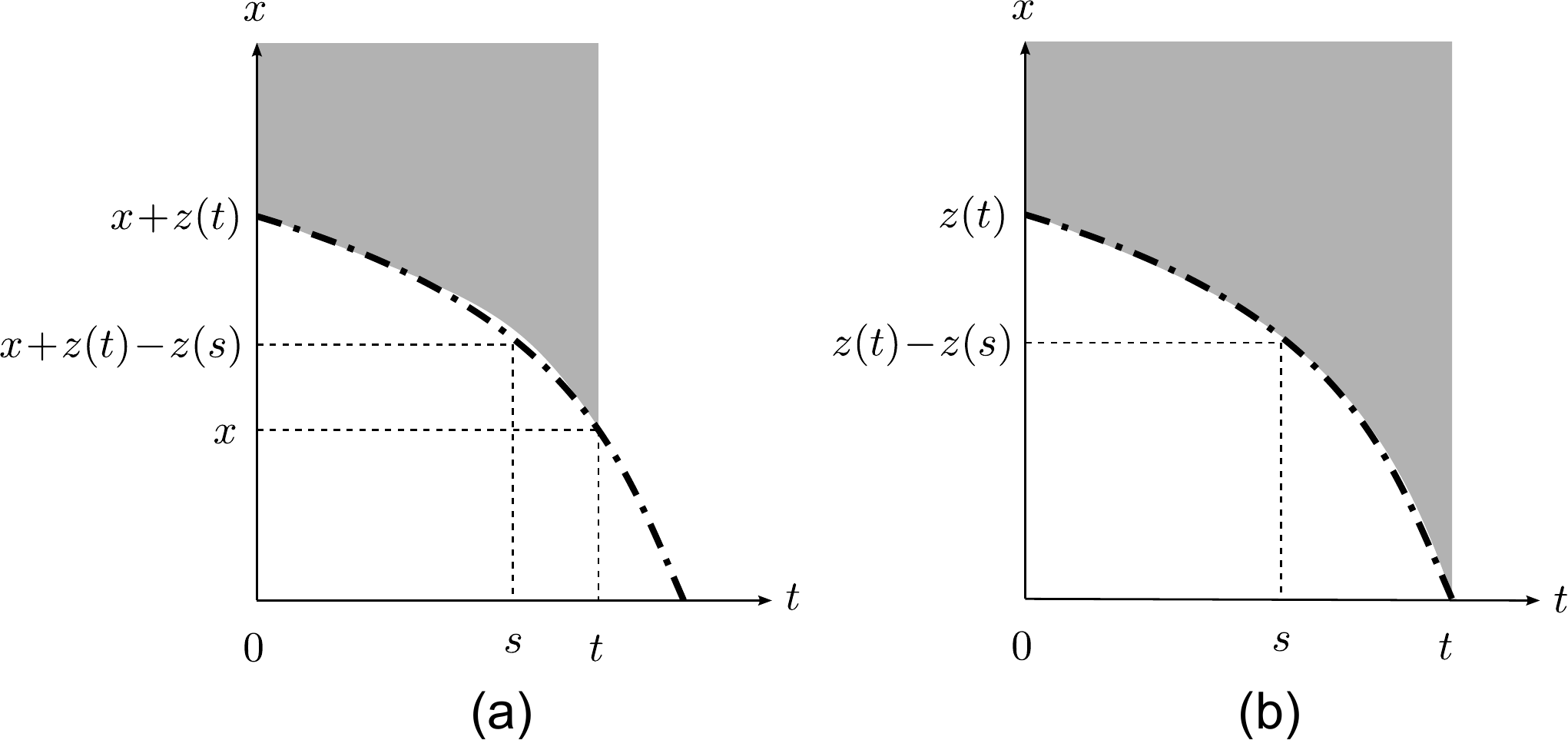}
 \caption{Illustration of the integral $K$-model } \label{integral-K-conservation}
 \ec\efg

In particular, for the initial value problem of \refe{K-model} with $f(t)=0$, \refe{K-model} can be simplified as
\bqn
\pd {}t K(t,x) -V\left(\frac {K(t,0)}L \right) \pd {}x K(t,x)&=& 0, \label{K-model-ivp}
\eqn
which is solved by $K(t,x)=K(0,x+z(t))$; i.e., $K(t,x)$ remains constant along the characteristic curve. In additon, $\la(t)=K(0,z(t))$. 
As a special case, if the initial distribution of the trip distance is negative exponential with $\Phi(0,x)=e^{-\frac x B}$, where $B$ is the average trip distance, we then have $K(0,x)=\la(0) e^{-\frac xB}$, $\la(t)=K(t,0)=\la(0) e^{-\frac {z(t)}B}$, and $K(t,x)=\la(t) e^{-\frac xB}$. Therefore, the remaining trip distance follows the same negative exponential distribution, and the average remaining trip distance is time-independent and still equals $B$.

\reff{K-model-ivp-sol} presents the analytical solution of \refe{K-model-ivp}, when the speed-density relation is piecewise constant as shown in \reff{K-model-ivp-sol}(a). The concave solid and dashed curves in \reff{K-model-ivp-sol}(b) are the characteristic curves. In addition, $K(t_1,0)=Lk_2$, $K(t_2,0)=Lk_1$, and $K(t_3,0)=0$; i.e., the completion times of trips $L k_2$, $L k_1$, and $0$ are respectively $t_1$, $t_2$, and $t_3$. Therefore, the travel speeds are respectively $v_3$, $v_2$, and $v_1$ during time intervals of $[0,t_1]$, $(t_1,t_2]$, and $(t_2, t_3]$, as shown in the figure. Along the characteristic curves, $x$ always decreases in $t$, since $v(t)>0$. As the speed is non-increasing in the density, the characteristic curves are concave. From the figure, we can also see that $K(t,x)$ decreases in both $t$ and $x$.  We can see that trips with different distances exit the network at different times, even though they all start at $t=0$. Therefore, the first-in-first-out principle does not apply in this model, but the shorter-(effective) distance-first-out principle still applies. Note that the violation of the first-in-first-out principle is not contradictory to other traffic flow models that satisfy the first-in-first-out principle on the same road, as the trips can have totally different origins and destinations. 

\bfg\bc
\includegraphics[width=5in]{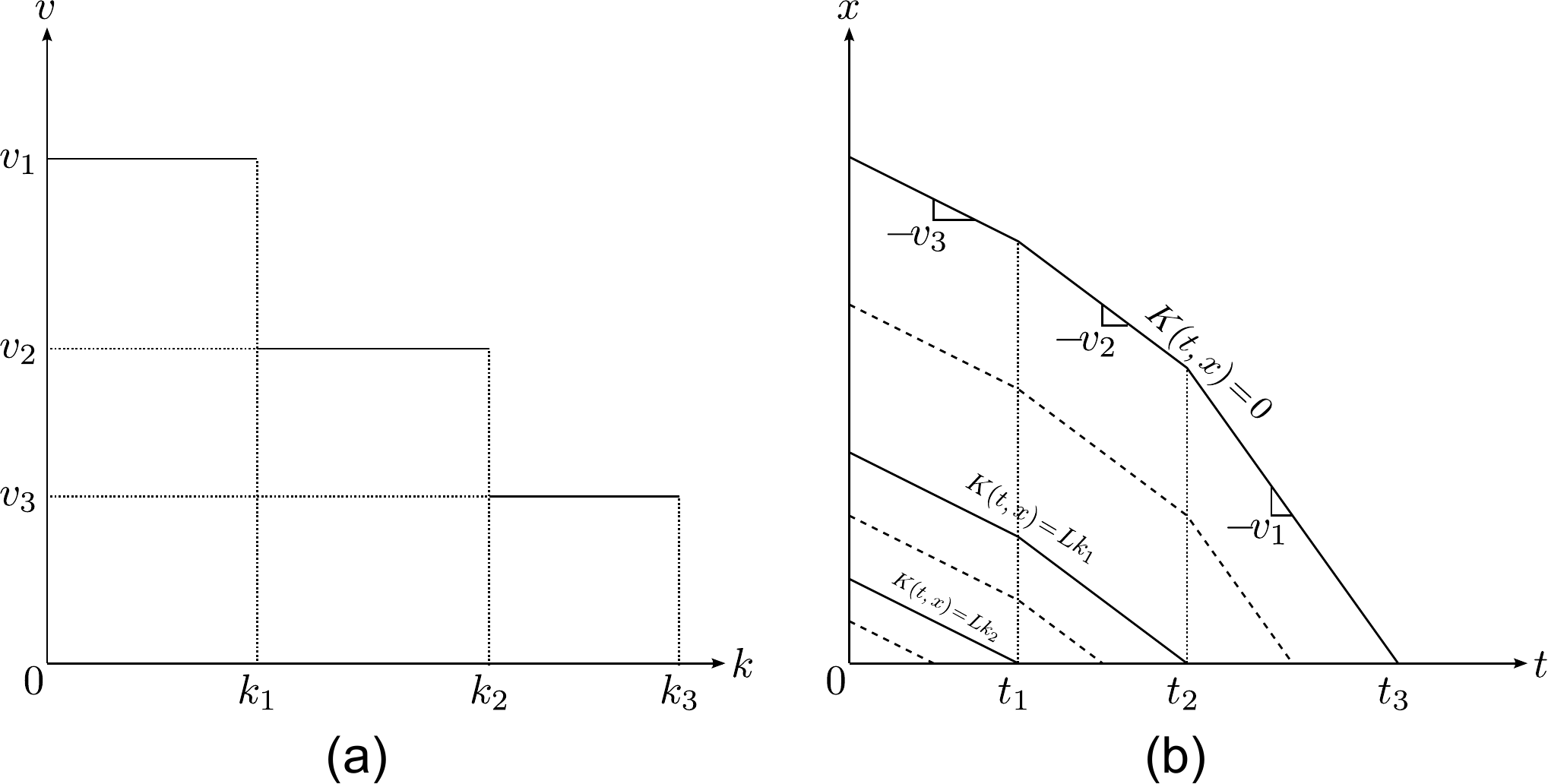}
\caption{Analytical solution of the initial value problem for the $K$-model}\label{K-model-ivp-sol}
\ec\efg

\subsection{The bathtub model in the $(z,x)$ coordinates}
When $v(t)>0$, $t$ strictly increases in $z(t)$. Thus the $(t,x)$ coordinates can be replaced by the $(z,x)$ coordinates. \footnote{We can also use $(\theta,x)$ coordinates for the bathtub model.}
In the new coordinates, the $K$-model, \refe{K-model}, can be re-written as
\bqn
\pd {}z K(z,x)-\pd{}x K(z,x)&=&\frac{f(z) \tilde \Phi(z,x)}{V\left(\frac{K(z,0)}L \right)}, \label{K-model-z}
\eqn
where $K(z,x)=K(\tau(z),x)$, $f(z)=f(\tau(z))$, and $\tilde \Phi(z,x)=\tilde \Phi(\tau(z),x)$. In the new coordinates, the characteristic curve is given by $\dot x(z)=-1$.
\refe{def:zt}, \refe{speed-density-relation},  and \refe{integral-lambda} can be re-written as 
\bsq \label{la-v-t}
\bqn
\la(z)&=& K(0,z)+\int_0^z \frac {f(y) \tilde \Phi(y, z-y)} { v(y)}  dy,\\
v(z)&=&V\left(\frac{\la(z)}L \right),\\
\tau(z)&=&\int_0^z \frac 1 { v(y)} dy,
\eqn
\esq
where $\la(z)=\la(\tau(z))$, and $v(z)=v(\tau(z))$.  \refe{integral-bathtub} can be re-written as
\bqn
K(z,x)&=&K(0,x+z)+\int_0^z \frac{f(y) \tilde \Phi(y, x+z-y)}{v(y)} dy. \label{integral-bathtub-z}
\eqn
Therefore, the number of trips with effective distances greater than $x+z$ at $z$, $K(z,x)$, equals the number of initial trips with distances larger than $x+z$, $K(0,x+z)$, plus those entering before $z$ with effective distances greater than $x+z$, $\int_0^z \frac{f(y) \tilde \Phi(y, x+z-y)}{v(y)} dy$.

Then 
\bqs
F(z)&=&F(\tau(z))=\int_0^z \frac{f(y)}{v(y)} dy, 
\eqs
or equivalently
\bqn
\dot F(z)&=&\frac{f(z)}{V\left(\frac{\la(z)}L \right)}, \label{Fz-ode}
\eqn
and 
\bqs
G(z)&=&\la(0)+F(z)-\la(z)=\la(0)-K(0,z)+\int_0^z \frac{f(y)}{v(y)}(1-\tilde \Phi(y, z-y)) dy.
\eqs
Thus $G(z)$ includes all initial trips whose distances are shorter than $z$, $\la(0)-K(0,z)$, and all entering trips whose effective distances are shorter than $z$, $\int_0^z \frac{f(y)}{v(y)}(1-\tilde \Phi(y, z-y)) dy$.

In the new coordinates, a trip can be represented by the cumulative travel distance when it enters the network, $z$, and its distance $x$. For such a trip, its effective distance is $x+z$, and its travel time is
\bqs
\Upsilon(z,x)&=&\Upsilon(\tau(z),x)=\tau(z+x)-\tau(z).
\eqs 
The average travel time for trips entering at $z$ is
\bqs
\bar \Upsilon(z)&=&\int_0^\infty \Upsilon(z,x) \tilde \varphi(z,x)dx=\int_0^\infty \tilde \Phi(z,x) \frac 1{v(z+x)} dx,
\eqs
where $\tilde \varphi(z,x)=\tilde \varphi(\tau(z),x)$.

\subsection{Finite difference methods and a numerical example}
The integral form of the bathtub model can be solved numerically. In particular, a difference-integration method can be developed to solve $\la(t)$ and $z(t)$ from \refe{integral-lambda}, and then  $K(t,x)$ and $G(t)$ from \refe{integral-bathtub} and  \refe{conservation-total-trips} respectively as follows. We divide the range of trip distances $[0,X]$ into $I$ intervals with $\dx= \frac X {I}$, where $X$ is the maximum trip distance.  We discretize the study period $[0,T]$ to $J$ time steps with $\dt=\frac TJ$. At $j\dt$ ($j=0,1,\cdots,J$),  the number of active trips in the network is $\la^j$, the travel speed by $v^j$, and the cumulative travel distance by $z^j$.
\bsq\label{intral-numerical-method}
\bqn
v^j &=& V\left( \frac {\la^j} L \right),\\
z^{j+1} &=& z^j +v^j \dt,\\
\la^{j+1} &=&K(0, z^{j+1}) + \sum_{i=0}^ j f(i\dt) \tilde \Phi(i\dt, z^{j+1} -z^i) \dt.
\eqn
\esq
Then for $j=0,\cdots, J$ and $i=0, \cdots, I$ 
\bqs
K(j\dt, i\dx)&=&K(0,i\dx+z^{j})+\sum_{m=0}^{j-1} f(m\dt) \tilde \Phi(m\dt, i\dx+ z^{j} -z^m ) \dt,\\
G(j\dt) &=& \la^0+ \sum_{m=0}^{j-1} f(m\dt) \dt -\la^{j}.
\eqs
A similar method can be developed in the $(z,x)$ coordinates.

Another method is to directly solve $K(t,x)$ from the differential form of the bathtub model, \refe{K-model}.
At $j$th time step, the time is denoted by $t_j$ ($j=0,1,\cdots$), where $t_0=0$. At $t_j$, $K(t_j,i\dx)$ is denoted by $K_i^j$ ($i=0,\cdots,I$),  and the number of active trips in the network is $\la^j=K_0^j$. Further the corresponding speed is denoted by $v^j=V(\frac{\la^j}L)$. If we denote the step-size $\dt^j= \frac \dx {v^j}$, then $t_{j+1}=t_j+\dt^j$, and $z(t_j)=j\dx$. From \refe{discrete-K-eqn} we can update $K_i^{j+1}$ ($i=0,\cdots,I-1$) as follows:
\bqn
K_i^{j+1}&=&K_{i+1}^j+\tilde K_i^j,
\eqn
where  $\tilde K_i^j=f(t_j) \tilde \Phi(t_j, i \dx) \dt^j$.\footnote{The finite difference method is effectively implemented in the $(z,x)$ coordinates. A related method is to first solve $\la(z)$, $v(z)$, and $\tau(z)$ by discretizing \refe{la-v-t} and then solve \refe{integral-bathtub-z}.}  The algorithm stops when the network is gridlocked with $v^j=0$ or when $t_{j+1}$ is beyond the study period $T$. 

As an example, we consider a network with $L=10$ lane-miles, the speed-density relation is \footnote{The network resembles the  road network of the Ipswich town center considered in \citep{godfrey1969mechanism}.}
\bqs
V(\r)&=&\min\{30, \frac{750} \r, 10(\frac{200} \r -1) \}.
\eqs
 Initially the network is empty,  the in-flux during a peak period is given by the following trapezoidal function:
\bqs
f(t)&=&\max\{0, \min\{ 10000 t, 4000, 10000(1-t) \}\},
\eqs
and the trip distance follows a uniform distribution:
\bqs
\tilde \Phi(t,x)&=&\max\{0, 1-\frac x{2 \tilde B(t)} \},
\eqs
where $\tilde B(t)$ is the time-dependent average trip distance at $t$:
\bqs
\tilde B(t)&=&2+\max\{0, \min\{ 7.5 t, 3, 7.5(1-t) \}\}.
\eqs
Here the maximum remaining trip distance is $X=5$ miles, and both the in-flux and the average trip distance are symmetric and reach their respective maximum values between 0.4 and 0.6 hr. The initial condition is given by $K_i^0=0$. We simulate the system until the accumulative travel distance reaches 30 miles.
We first let $\dx=2^{-6}$ mile and $I=320$. The evolution of $K(t,x)$ is illustrated in \reff{generalized_bathtub_example}(a), which shows that the number of active trips reaches the maximum value between 0.75 and 1 hr. Thus the network is most congested between 0.75 and 1 hr, even though the demand pattern peaks between 0.4 and 0.6 hr. In \reff{generalized_bathtub_example}(b) we present the solutions of $z(t)$ for different $\dx$'s. From this figure we can see that the solutions converge with diminishing $\dx$; it suggests that the numerical solution converges to the theoretical one when $\dx$ diminishes. By comparing the times for the cumulative travel distances to reach 30 miles, we obtain an approximate convergence rate of 1.

\bfg\bc
$\ba{c@{\hspace{0.3in}}c}
\includegraphics[height=2in]{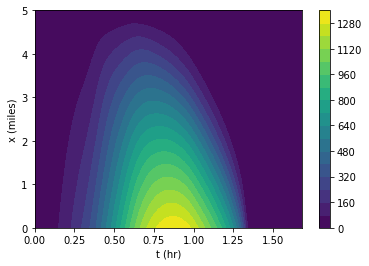} &
\includegraphics[height=2in]{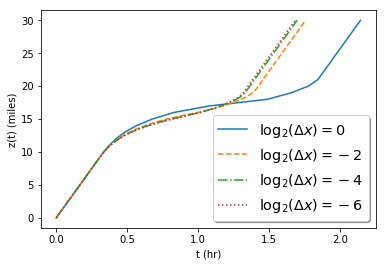} \\
\mbox{ (a) Solution of $K(t,x)$ with $\dx=2^{-6}$ mile} &
\mbox{ (b) Solutions of $z(t)$ for different $\dx$'s}
\ea$
\caption{A numerical example of the generalized bathtub model} \label{generalized_bathtub_example}
\ec\efg

\section{Special case 1: Vickrey's bathtub model}
 \citet{vickrey1991congestion,vickrey1994types} presented a special case of the bathtub model when the entering trips' distances and the remaining trip distance follow the same time-independent negative exponential distribution 
\bqn
\Phi(t,x)&=&\tilde \Phi(t,x)=e^{-\frac xB}, \label{negative-exponential}
\eqn
where $B$ is the average trip distance and the average remaining trip distance.

With the negative exponential distribution of trip distances, the only unknown variable is $\la(t)$, and the bathtub model is equivalent to the simplified version of \refe{la-ode}:
\bqn
\dot \la (t)&=& f(t) - \frac 1B \la(t) V\left(\frac{\la(t)}L \right). \label{Vickrey-bathtub}
\eqn
 In this case, the out-flux or the trip completion rate equals  \citep{gonzales2012morning}
\bqn
g(t)&=&\frac 1B \la(t) v(t)=\frac LB Q(\r(t)), \label{vickrey-outflux}
\eqn
which is proportional to the traffic flow rate  $Q(\r(t))$.
We can see that Vickrey's bathtub model for trips served by privately operated vehicles is derived based on three premises: (i) the bathtub treatment of a road network, (ii) the speed-density relation at the network level, and (iii) the time-independent negative exponential distribution of trip distances. The generalized bathtub model in \refe{K-model} relaxes the third premise. In addition, Vickrey's bathtub model tracks the evolution of the total number of trips, but the generalized bathtub model tracks the evolution of the number of trips with a remaining distance over any value.

In \citep{arnott2016equilibrium,arnott2018solving}, Vickrey's bathtub model was considered improper as ``the exit rate from downtown traffic (i.e., the arrival rate at work) depends only on the density of downtown traffic, and hence that the first exit occurs as soon as the first entry''.  Even though these phenomena may not occur when the trip distances follow other distributions, they are physically meaningful: (i) since the remaining trip distance distribution is time-independent, the out-flux only depends on the density; and (ii) since some trips have zero distances with the negative exponential distribution, these trips can exit the network in no time.

\subsection{Properties}
The properties of the stationary states, including the gridlock state, and their stability properties in Section \ref{section:stationarystate} can be easily extended for Vickrey's bathtub model. 
Stability and control of \refe{Vickrey-bathtub} were discussed in \citep{agnew1976dynamic,daganzo2007gridlock}.

From \refe{integral-lambda}, we have the following integral form of Vickrey's bathtub model
\bqn
\la(t)&=&\la(0) e^{-\frac{z(t)}B}+\int_0^t f(s)e^{-\frac{z(t)-z(s)}B}  ds, \label{integral-Vickrey}
\eqn
which and $\dot z(t)=V(\frac{\la(t)}L)$ can be used to solve $\la(t)$ and $z(t)$,
and
\bqn
G(t)&=&(1-e^{-\frac{z(t)}B})\la(0)+\int_0^t f(s)(1-e^{-\frac{z(t)-z(s)}B}) ds. \label{G-Vickrey}
\eqn
In the $(z,x)$ coordinates, Vickrey's bathtub model can be written as
\bqn
\dot \la(z)&=&\frac{f(z)}{v(z)}-\frac 1B \la(z),
\eqn
whose integral form is
\bqn
\la(z)&=&e^{-\frac zB} [\la(0)+\int_0^z \frac{f(y) } {v(y)} e^{\frac yB}dy ].
\eqn
Thus the cumulative out-flow is
\bqn
G(z)&=&(1-e^{-\frac zB})\la(0)+\int_0^z \frac{f(y) } {v(y)} (1-e^{\frac {y-z}B})dy.
\eqn
With $\la(z)$ we can calculate $\tau(z)$ from \refe{def:tau} and the trip travel times from \refe{def:traveltime}.
Under general initial and boundary conditions, Vickrey's bathtub model and its integral form in both $t$ and $z$ coordinates can be numerically solved with, for example, Euler's forward method. It is easy to show that \refe{intral-numerical-method} is equivalent to Euler's forward method for \refe{Vickrey-bathtub}; i.e., the numerical methods for both the integral and differential forms are equivalent for Vickrey's bathtub model. The computational cost is proportional to $\frac 1\dt$.

As shown in Section \ref{section:general:analytical}, the initial value problem for Vickrey's bathtub model, \refe{Vickrey-bathtub},  with $\la(0)>0$ can be solved by 
\bqn
\la(t)&=&\la(0) e^{-\frac{z(t)}B},
\eqn
and the average travel time for the initial trips is given by
\bqn
\bar \Upsilon(0)&=&\int_0^\infty \frac{e^{-\frac xB} }{V\left( \frac{\la(0)}L e^{-\frac xB}\right)} dx.
\eqn
In this case, \refe{Vickrey-bathtub} is equivalent to
\bqn
\dot z(t) &=&V\left( \frac{\la(0)}L e^{-\frac {z(t)}B}\right),
\eqn
which can be analytically solved with a given speed-density relation.

Next we solve Vickrey's bathtub model with a constant in-flux on an initially empty road network. That is, $f(t)=f$, and $\la(0)=0$. When the demand level is relatively low; i.e., when \refe{gridlock-sufficientcondition} is violated, 
$\la(t)$ keeps increasing until it reaches a stationary $\la$ such that $Bf= L Q\left( \frac\la L\right)$. However, if \refe{gridlock-sufficientcondition} is satisfied; i.e., if the demand level is too high, the network becomes gridlocked eventually.

\subsection{Equivalent conditions for Vickrey's bathtub model}

\begin{theorem} \label{thm:vickreymodel}
	For a network, which is never gridlocked with $v(t)>0$ for any $t$, the following six statements are equivalent and lead to Vickrey's bathtub model:
	\bi
	\item [(i)] The entering trips' distances and the remaining trip distance follow the same time-independent negative exponential distribution as in \refe{negative-exponential}. 
	
	\item [(ii)] The remaining trip distance follows a time-independent negative exponential
	\bqn
	\Phi(t,x)&=&e^{-\frac xB}.
	\eqn
	\item [(iii)] The distributions of the initial trips' remaining distances and the entering trips' total distances follow the same, time-independent negative exponential distribution:
	\bqn
	\Phi(0,x)&=&\tilde \Phi(t,x)=e^{-\frac xB}.
	\eqn
	\item [(iv)] $\Phi(t,x)=\tilde \Phi(t,x)=\Phi(x)$ are time-independent and continuous.
	\item [(v)] $B(t)=\tilde B(t)=B$ is time-independent.
	\item [(vi)] The out-flux can be written as $g(t)=\frac 1{B} \la(t) v(t)$; i.e., $\varphi(t,0)=\frac 1B$ is constant.
	\ei
	In particular, the statements are equivalent in the following two special cases: (i) the initial network is empty; i.e., $\la(0)=0$; (ii) no trips enter the network at any time; i.e., $f(t)=0$.
\end{theorem}
{\em Proof}. Clearly (i) $\Rightarrow$ (ii), (iii), (iv), (v), and (vi). In addition, (iv) $\Rightarrow$ (v) from the definitions of $B(t)$ and $\tilde B(t)$.

If $\Phi(t,x)=e^{-\frac xB}$, from \refe{K-model} we have
\bqs
\dot \la(t) e^{-\frac xB}+\frac {\la(t)v(t)} B e^{-\frac xB} &=&f(t) \tilde \Phi(t,x).
\eqs
Letting $x=0$, we have
\bqs
\dot \la(t) +\frac {\la(t) v(t)} B  &=&f(t) .
\eqs
Comparing the above two equations, we have $\tilde \Phi(t,x)=e^{-\frac xB}$. Therefore, (ii) $\Rightarrow$ (i) and (iii).

We prove that (iii) $\Rightarrow$ (ii) by induction with the discrete $\Phi$-model in \refe{discrete-cdf-eqn}.  First, $\Phi(0,x)=e^{-\frac xB}$ is given. Second, we assume that $\Phi(t,x)=\tilde \Phi(t,x)=e^{-\frac xB}$ at $t$. We then prove that $\Phi(t+\dt,x)=e^{-\frac xB}$ for an arbitrarily small $\dt$. From \refe{discrete-cdf-eqn} we have
\bqs
\la(t+\dt)\Phi(t+\dt,x)&=&\la(t) e^{-\frac{x+v(t)\dt}B}+f(t)\dt e^{-\frac xB}.
\eqs
Further setting $x=0$ in \refe{discrete-cdf-eqn}, we have
\bqs
\la(t+\dt)&=&\la(t) e^{-\frac{v(t)\dt}B}+f(t)\dt.
\eqs
Comparing the two equations, we can see that $\Phi(t+\dt,x)=e^{-\frac xB}$.
By induction, we conclude that $\Phi(t,x)=e^{-\frac xB}$ at any time $t$. Hence (iii) $\Rightarrow$ (ii).

Substituting  $\Phi(t+\dt,x)=\Phi(t,x)=\tilde \Phi(t,x)=\Phi(x)$ into \refe{discrete-cdf-eqn}, we have
\bqs
\la(t+\dt)\Phi(x)&=&\la(t) \Phi(x+v\dt)+f(t)\dt \Phi(x).
\eqs
Letting $x=0$, the above equation yields
\bqs
\la(t+\dt)&=&\la(t) \Phi(v\dt)+f(t)\dt.
\eqs
Comparing the above two equations we have
\bqs
\Phi(x+v\dt)&=&\Phi(x)\Phi(v\dt).
\eqs
Assuming $\Phi(1)=e^{-\frac 1B}$,  from the above equation we have
\bqs
\Phi(1)&=&\Phi(v\dt)^{\frac 1{v\dt}},
\eqs
which leads to $\Phi(t,v\dt)=e^{-\frac{v\dt}B}$ or $\Phi(t,x)=e^{-\frac xB}$. Therefore (iv) $\Rightarrow$ (ii).

If $B(t)=\tilde B(t)=B$, then from the definitions of $B(t)$ and $\tilde B(t)$ we can see that both $\Phi(t,x)=\Phi(x)$ and $\tilde \Phi(t,x)=\tilde \Phi(x)$ are time-independent. From \refe{K-model} we have
\bqs
\dot \la(t) \Phi(x)-v(t) \la(t) \dot \Phi(x)&=&f(t) \tilde \Phi(x).
\eqs
Integrating both sides of the above equation for $x$ from $0$ to $\infty$, we have 
\bqs
B \dot \la(t)+v(t) \la(t)&=&f(t) B.
\eqs
Comparing the above equation with \refe{continuous-la-eqn}, we can see that $g(t)=\frac 1B \la(t) v(t)$ or, equivalently, $\varphi(t,0)=\frac 1B$. Thus (v) $\Rightarrow$ (vi).

When $g(t)=\frac 1B \la(t)v(t)$, \refe{continuous-la-eqn} can be written as
\bqs
\dot \la(t)+ \frac 1B v(t) \la(t)&=&f(t),
\eqs
whose integral form is
\bqs
\la(t)&=&\la(0) e^{-\frac{z(t)}B}+\int_0^t f(s)e^{-\frac{z(t)-z(s)}B}  ds.
\eqs
Comparing the above equation with \refe{integral-lambda}, we can see that $\Phi(t,x)=\tilde \Phi(t,x)=e^{-\frac xB}$. Thus (vi) $\Rightarrow$ (i).

Therefore the equivalence among the six statements is proved. \eop

That (v) implies (ii) was first proved in \citep{vickrey1991congestion} under the condition that the in-flux $f(t)=0$ by observing that $K(t,x) dx=-B dK(t,x)$.

Note that (v) is not equivalent to that $\tilde B(t)=B$ is time-independent. A counter example is when all trip distances are equal as discussed in Section \ref{section:constanttripdistance}.

Also note that the equivalence among the six statements is independent of the speed-density relation as long as $v(t)$ is the same for all trips at the same time.

\section{Special case 2: Deterministic trip distances}
If the distance for trips entering the network at $t$ is deterministic and equals $\tilde B(t)$, then the proportion of the entering trips with distances not smaller than $x$ is
\bqn
\tilde \Phi(t,x)&=&H(\tilde B(t)-x)=\cas{{ll} 0, & x> \tilde B(t); \\ 1, & x\leq \tilde B(t),}
\eqn
where $H(\cdot)$ is the Heaviside function, and 
the $K$-model can be written as
\bqn
\pd {}t K(t,x) -V\left(\frac {K(t,0)}L \right) \pd {}x K(t,x)&=& f(t) H(\tilde B(t)-x), \label{K-model-deterministic}
\eqn
The integral $K$-model can be written as
\bqn
K(t,x)&=&K(0,x+z(t))+\int_0^t f(s) H(\tilde B(s)+z(s)-x-z(t)) ds, \label{integral-bathtub-deterministic}
\eqn
The integral form of $\la(t)$ is
\bqn
\la(t)&=&K(0,z(t))+\int_0^t f(s) H(\tilde B(s)+z(s)-z(t)) ds. \label{integral-la-deterministic}
\eqn
In this case, from \refe{def:triptraveltime} and \refe{def:traveltime1} we have $\Upsilon(t,x)=\Upsilon(t,\tilde B(t))=\bar \Upsilon(t)$. 

In the $(z,x)$ coordinates, the $K$-model, \refe{K-model}, can be re-written as
\bqn
\pd {}z K(z,x)-\pd{}x K(z,x)&=&\frac{f(z) H(\tilde B(z)-x)}{V\left(\frac{K(z,0)}L \right)}, \label{K-model-z-de}
\eqn
where $\tilde B(z)=\tilde B(\tau(z))$. The corresponding integral form is given by
\bqn
K(z,x)&=&K(0,x+z)+\int_0^z \frac{f(y) H(\theta(y)-x-z)}{v(y)} dy, \label{integral-bathtub-z-de}
\eqn
where $\theta(z)=\tilde B(z)+z$ is the effective distance of trips entering at $z$. 
Furthermore $v(z)$ is solved by
\bsq \label{la-v-t-de}
\bqn
\la(z)&=& K(0,z)+\int_0^z \frac {f(y) H(\theta(y)-z)} { v(y)}  dy,\\
v(z)&=&V\left(\frac{\la(z)}L \right),\\
\tau(z)&=&\int_0^z \frac 1 { v(y)} dy.
\eqn
\esq
In addition,  the cumulative out-flow at $z$ is
\bqn
G(z)&=&\la(0)-K(0,z)+\int_0^z \frac{f(y) H(z-\theta(y))}{v(y)} dy.
\eqn

\subsection{When $\der{}z \tilde B(z) <,=,> -1$}
When $\der{}z \tilde B(z) = -1$, $\der{}t \tilde B(t) = -v(t)$. In this case, $\theta(z)=\tilde B(0)$, and 
\bsq\label{la-deterministic}
\bqn
\la(z)&=& K(0,z)+ H(\tilde B(0)-z) \int_0^z \frac {f(y) } { v(y)}  dy=K(0,z)+ H(\tilde B(0)-z)  F(z).
\eqn 
 In addition, $K(z,x)=K(0,x+z)+H(\tilde B(0)-x-z) F(z)$. The cumulative out-flow at $z$ is $G(z)=\la(0)-K(0,z)+H(z-\tilde B(0))F(z)$.
In this case, all entering trips exit at the same time corresponding to the cumulative travel distance of $z=\tilde B(0)$.

When $\der{}z \tilde B(z)<-1$, $\der{}t \tilde B(t)<-v(t)$. In this case, $\der{}z \theta(z)<0$; i.e., the entering trips' effective distances decrease with $z$, and the trips follow the last-in-first-out principle. $\tilde B^{-1}(0)$ is the cumulative travel distance when the entering trip's distance becomes zero. The bathtub model has the following three types of solutions. (i) When $z<\tilde B^{-1}(0)$, $\theta(y)> \theta(\tilde B^{-1}(0))=\tilde B^{-1}(0)$ for $y\in[0,z]$. Thus, no entering trips exit the network; i.e., $G(z)=\la(0)-K(0,z)$, and $\la(z)=K(0,z)+F(z)$. (ii) When $\tilde B^{-1}(0)\leq z \leq \tilde B(0)$, $G(z)=\la(0)-K(0,z)+F(\tilde B^{-1}(0))-F(\theta^{-1}(z))$, and $\la(z)=K(0,z)+F(\theta^{-1}(z))$. (iii) When $z>\tilde B(0)$, all entering trips exit the network, and $G(z)=\la(0)-K(0,z)+F(z)$, and $\la(z)=K(0,z)$.
Thus $\la(z)$ is solved by
\bqn
\la(z)=\cas{{ll} K(0,z)+F(z), & z<\tilde B^{-1}(0); \\ K(0,z)+F(\theta^{-1}(z)), & \tilde B^{-1}(0)\leq z \leq \tilde B(0); \\ K(0,z), & z>\tilde B(0).}
\eqn

When $\der{}z \tilde B(z) >-1$, $\der{}t \tilde B(t)>-v(t)$. In this case, $\dot \theta(z)>0$, and $\theta(z)\geq\theta(0)=\tilde B(0)$ for $z\geq0$; i.e., the entering trips' effective distances increase with $z$, and the trips follow the first-in-first-out principle. In this case, there are two types of solutions. (i) When $z<\tilde B(0)$, $\theta(y)>z$ for $y\in[0,z]$. Thus $G(z)=\la(0)-K(0,z)$; i.e., only initial trips exit the network, since entering trips' effective distances are greater than $z$. In addition, $\la(z)=K(0,z)+F(z)$. (ii) When $\tilde B(0) \leq z \leq  \tilde B^{-1}(0) $, all trips entering before $\theta^{-1}(z)$ exit the network. (iii) When $\tilde B^{-1}(0)<\infty$ and $z>\tilde B^{-1}(0)$, $\la(z)=K(0,z)$. Thus $G(z)=\la(0)-K(0,z)+F(\theta^{-1}(z))$, and $\la(z)=K(0,z)+F(z)-F(\theta^{-1}(z))$.
Thus $\la(z)$ is solved by
\bqn
\la(z)&=&\cas{{ll} K(0,z)+F(z), & z<\tilde B(0); \\ K(0,z)+F(z)-F(\theta^{-1}(z)), & \tilde B(0) \leq z \leq  \tilde B^{-1}(0); \\ K(0,z), & z>\tilde B^{-1}(0).} \label{la-deterministic3}
\eqn
\esq

\refe{Fz-ode} and \refe{la-deterministic} form a deterministic bathtub model in the three cases. With $\la(z)$ we can calculate $\tau(z)$ from \refe{def:tau} and the trip travel times from \refe{def:traveltime}.
As a special case, when the entering trips' distances are constant $\tilde B(z)=\tilde B$.

The three cases are illustrated in \reff{bathtub_deterministic_sol}, where the dashed lines represent the characteristic curves in the $(z,x)$-space, and the solid lines for $x=\tilde B(z)$. In \reff{bathtub_deterministic_sol}(b) and (c), we also illustrate $\theta^{-1}(z)$. Note that, if we let $y=\theta^{-1}(z)$, then $z=\theta(y)$, which is the effective distance for the trip entering at $y$. From the figures, we can easily verify the solutions in the three cases. Also note that $\der{}z \tilde B(z)<-1$ in the second case corresponds to the regular sorting condition in \citep{fosgerau2015congestion}.

\bfg\bc
\includegraphics[width=6in]{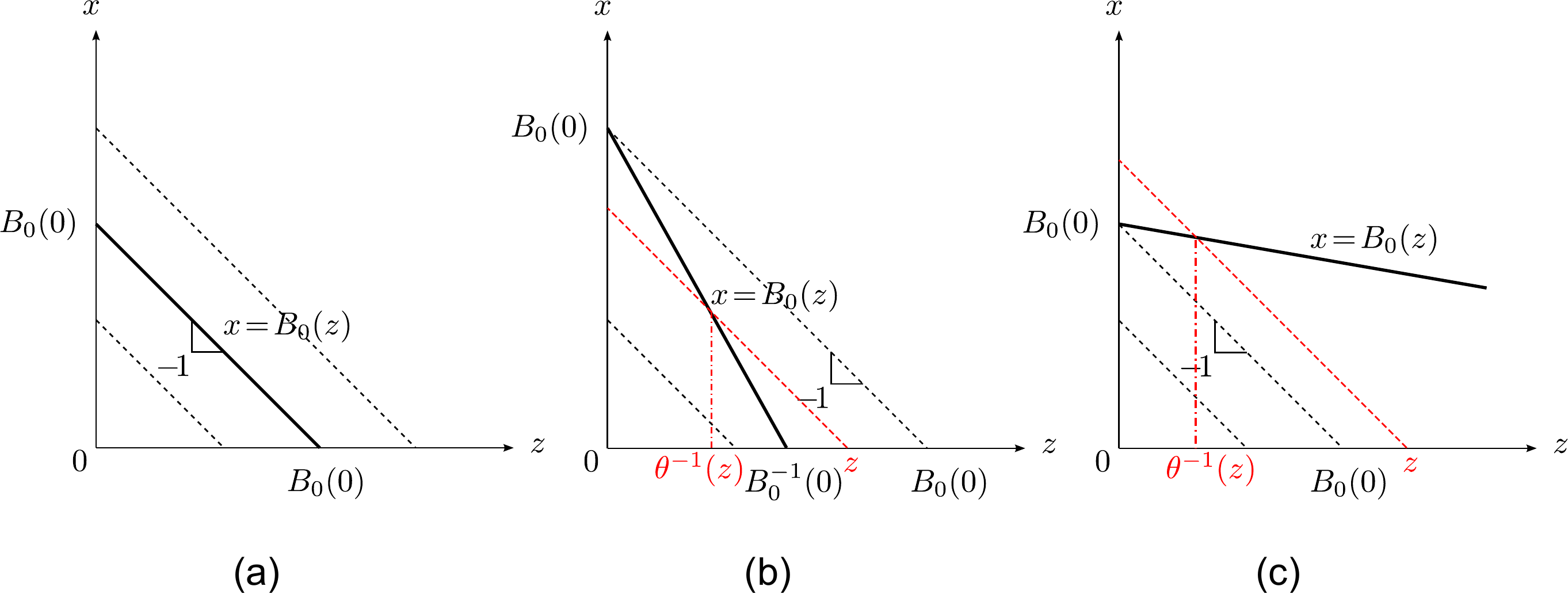}
\caption{Illustration of three cases of the bathtub model with deterministic trip distances: (a) $\der{}t \tilde B(t) = -v(t)$; (b) $\der{}t \tilde B(t) < -v(t)$; (c) $\der{}t \tilde B(t) > -v(t)$. }\label{bathtub_deterministic_sol}
\ec\efg

\subsection{Constant trip distances} \label{section:constanttripdistance}
In this subsection, we study the generalized bathtub model with constant trip distances $\tilde B(z)=\tilde B$ with initial trips' distances not greater than $\tilde B$. In this case, \refe{K-model-deterministic} can be simplified as ($x\in[0,\tilde B]$)
\bqn
\pd {}t K(t,x) -V\left(\frac {K(t,0)}L \right) \pd {}x K(t,x)&=& f(t) , \label{K-model-constant}
\eqn 
which cannot be simplified into Vickrey's bathtub model, \refe{Vickrey-bathtub}, since $\varphi(t,0)$ may not equal $\frac 1{\tilde B}$ at all times. As an example, for an initially empty network, $\varphi(0,0)=0$, and no vehicles can leave the network at $t=0$; but it is not the case in Vickrey's model, as those trips of zero distances can exit the network at $t=0$.

For $x\in[0,\tilde B]$ and $t\geq 0$, we define $N(t,x)$ as the cumulative number of trips passing $x$ by $t$. Since $\der{}t \tilde B(t)=0 >-v(t)$, trips follow the first-in-first-out principle,  $N(t,x)$ equals the number of trips which exit the network before the one passing $(t,x)$. We assume that the initial trips' distances are not greater than $\tilde B$. Thus 
\bqn
N(t,x)&=&N(0,\tilde B)+F(t)-K(t,x),
\eqn
and \refe{K-model-constant} is equivalent to
\bqn
\pd{}t N(t,x) - V\left(\frac{N(t,\tilde B)-N(t,0)}L \right) \pd {}x N(t,x)&=&0. \label{nonlocal-hj}
\eqn
From the definition of $N(t,x)$ we can see that $k(t,x)=\pd{}x N(t,x)$, and
\bqs
\pd{}t k(t,x) - V\left(\frac1L \int_0^{\tilde B} k(t,y) dy \right) \pd {}x k(t,x)&=&0, \label{nonlocal-lwr}
\eqs
which is the simplified version of the $k$-model, \refe{k-model}. 
For this model, the initial conditions are given by $k(0,x)$ or $N(0,x)$, and the boundary conditions are given by $k(t,\tilde B)=\frac{f(t)}{v(t)}$ or $N(t,\tilde B)=N(0,\tilde B)+F(t)$.
If the transportation network is mapped into a virtual road, whose length is $\tilde B$ and an average number of lanes of $\frac L {\tilde B}$, then $k(t,x)$ is the density of trips at time $t$ and location $x$, and $\frac1L \int_0^{\tilde B} k(t,y) dy$ is the average per-lane density on the road.   In this sense, \refe{nonlocal-lwr} is a nonlocal conservation law, which is a counterpart of the traditional LWR model \citep{lighthill1955lwr,richards1956lwr}, and \refe{nonlocal-hj} is a nonlocal Hamilton-Jacobi equation \citep{newell1993sim}. As the LWR model, \refe{nonlocal-lwr} also satisfies the conservation law, but the travel speed is not determined by the local density as in the LWR model; thus it is a nonlocal model.
In addition, the number of lanes on the virtual road may be inhomogeneous with respect to location and unbounded. But different from the LWR model, only the average number of lanes is relevant in the bathtub model.

\subsubsection{Other equivalent formulations}

We denote the traveling time of a trip entering the system at $t$ by $\bar \Upsilon(t)$. Thus
\bsq \label{cumulativeflow-model}
\bqn
N(0,\tilde B)+F(t)&=&G(t+\bar \Upsilon(t)),
\eqn 
and
\bqn
\int_t^{t+\bar \Upsilon(t)} v(s) ds&=&\tilde B.
\eqn
In addition, the number of vehicles in the system at $t$ equals $N(0,\tilde B)+F(t)-G(t)=N(t,\tilde B)-N(t,0)$, and
\bqn
v(t)&=&V\left(\frac{N(0,\tilde B)+F(t)-G(t)}{L}\right).
\eqn
\esq
Note that $N(0,\tilde B)$ and $F(t)$ are given. Thus \refe{cumulativeflow-model} is a new formulation of the bathtub model with constant trip distances, in which the unknown variables are $G(t)$, $\bar \Upsilon(t)$, and $v(t)$. It was derived as the ``proper bathtub model'' and solved as delay-differential equations in \citep{arnott2018solving}.

We denote the time for trip $n$ to pass $x$ on the virtual road by $T(n,x)$. Thus $T(n,x)=\tau(\tilde B-x+z(T(n,\tilde B)))$, which is also the time for trip $n$'s remaining distance to be $x$. In the $(n,x)$ coordinates, the bathtub model can be written as
\bqs
\int_{T(n,\tilde B)}^{T(n,x)} V\left(\frac{N(0,\tilde B)+F(s)-G(s)}{L}\right) ds =\tilde B-x,
\eqs
whose differential form is
\bqn
T_x(n,x)&=& - \frac 1{  V\left(\frac{N(0,\tilde B)+F(T(n,x))-G(T(n,x))}{L}\right)}.
\eqn
In this formulation, $T(n,x)$ is the unknown variable, and $G(t)=n$ if $T(n,0)=t$.

We denote trip $n$'s location at $t$ on the virtual road  by $X(t,n)$. Thus $X(t,n)=\tilde B-z(t)+z(T(n,\tilde B))$, which is the remaining distance of trip $n$ at $t$. In the $(t,n)$ coordinates, the bathtub model can be written as:
\bsq
\bqn
X_t(t,n)&=&-V\left(\frac{N(0,\tilde B)+F(t)-G(t)}{L}\right), \label{trajectory-model}
\eqn
where $X_t(t,n)$ is the travel speed, and
\bqn
G(t)&=&\max_n \{ n | X(t,n) \leq 0 \}.
\eqn
\esq
An integral formulation of \refe{trajectory-model} is
\bqs
\int_{T(n,\tilde B)}^t V\left(\frac{N(0,\tilde B)+F(s)-G(s)}{L}\right) ds =\tilde B-X(t,n),
\eqs
where $T(n,\tilde B)$ is the time for trip $n$ to enter the virtual road.
A special case of the above equation is 
\bqs
\int_{T(n,\tilde B)}^{T(n,0)} V\left(\frac{N(0,\tilde B)+F(s)-G(s)}{L}\right) ds =\tilde B,
\eqs
which is equivalent to (\ref{cumulativeflow-model}b).

\subsubsection{Properties and solutions}
In the stationary states, $\pd {} t K(t,x)=0$. Thus $\pd{}t k(t,x)=0$. From \refe{nonlocal-lwr}, we find that $k(t,x)=k$ is constant for non-gridlocked states. Hence the number of active trips $\la=\tilde B k$. Further from \refe{K-model-constant} we have $f(t)=\la V(\frac \la L)$. In this case, the active trip's remaining distance follows a uniform distribution: $\varphi(t,x)=\frac 1{\tilde B}$ for $x\in[0,\tilde B]$. 
The out-flux is $g(t)=\frac 1 {\tilde B} \la(t) v(t)$, which is the same as that in Vickrey's bathtub model, \refe{vickrey-outflux}; thus Vickrey's bathtub model provides a good approximation ``in steady state, or when the accumulation (number of active trips) varies slowly with slow demand variations'' \citep{mariotte2017macroscopic}.
This is a special case of those discussed in Theorem \ref{thm:stasta}. The stability property of the stationary states is also discussed in Section \ref{section:stationarystate}.

For an initially empty network with constant trip distances, \refe{la-deterministic3} can be simplified as
\bqn
\la(z)&=&\cas{{ll} F(z), & z<\tilde B; \\ F(z)-F(z-\tilde B), & z\geq \tilde B.}
\eqn
Then numerically we can solve $\la(z)$, $F(z)$, and $\tau(z)$ as follows. First, we divided $[0,\tilde B]$ into $I$ intervals with $\dz=\frac{\tilde B}I$. We denote $\la_j=\la(j\dz)$, $F_j=F(j\dz)$, $f_j=f(j\dz)$, and $\tau_j=\tau(j\dz)$. Given the initial condition, $\la_0=0$ and $F_0=0$, and the boundary condition, $f_j$, we can then calculate $F_{j+1}$, $\la_{j+1}$, and $\tau_{j+1}$ as follows
\bsq
\bqn
\Delta \tau_j&=&\frac{\dz}{V(\frac{\la_j}L)},\\
\tau_{j+1}&=&\tau_j+\Delta \tau_j,\\
F_{j+1}&=&F_j+\Delta \tau_j f_j,\\
\la_{j+1}&=&\cas{{ll} F_{j+1}, & j<I-1; \\ F_{j+1}-F_{j+1-I}, & j\geq I-1,}
\eqn
\esq
where $\Delta \tau_j$ is the time for the cumulative travel distance to increase from $j\dz$ to $(j+1)\dz$. From $z(t)$ and $\tau(t)$ we can then calculate $\bar \Upsilon(t)=t-\tau(z(t)-\tilde B)$,  and $N(t,x)=F(\tau(x+z(t)-\tilde B))$ for any $x\in[0,\tilde B]$ and $t\geq 0$. Similarly we can also calculate $X(t,n)$ and $T(n,x)$. An equivalent numerical method in the $(t,x)$ coordinates can be simplified from \refe{intral-numerical-method} but omitted here.

\section{Extensions for other types of trips}
In the preceding sections, the conservation equations of $\la(t)$ and $K(t,x)$ in \refe{continuous-la-eqn} and \refe{continuous-K-eqn} respectively, as well as their respective integral forms (\ref{integral-lambda}a) and \refe{integral-bathtub}, are coupled with the speed-density or speed-trip relation in \refe{speed-density-relation} to describe the evolution of trips served by privately operated vehicles.
We can see that, as long as all trips have the same travel speed, $v(t)$, different formulations of the generalized bathtub model in \refe{continuous-la-eqn}, \refe{continuous-K-eqn}, (\ref{integral-lambda}a), and \refe{integral-bathtub} are still valid. Therefore, the bathtub model and its properties and solutions can be directly extended when $v(t)$ is exogenously given or depends on other variables.

\subsection{Trips served by mobility service vehicles}
In the generalized bathtub model, \refe{K-model}, the vehicles and trips are integrated, and the speed-density relation in \refe{speed-density-relation} directly determines the speed from the number of active trips. But when trips are served by mobility service vehicles, the vehicles and trips are disintegrated, and \refe{la2rho} and \refe{speed-density-relation} may not be true. That is, vehicle density per lane $\rho(t)$ may be independent of the number of active trips due to varying occupancies of vehicles. Therefore we extend the  fundamental diagram of network trip flows when trips are served by mobility service vehicles as
\bqn
v(t)&=&V(\rho(t), \la(t), f(t), g(t)). \label{nfd-mobility-service}
\eqn
Clearly the number of running vehicles on the network could impact the travel speed due to the vehicle congestion effect. In a dedicated mobility service system, such as a metro system or a bus rapid transit system, $\rho(t)$ only includes the mobility service vehicles; but in a normal bus network or a mobility service system with transportation network company vehicles, $\rho(t)$ also includes other privately operated vehicles.
In addition, the number of passengers could also impact the average travel speed, due to different boarding and alighting times, which can depend on $\la(t)$, $f(t)$, and $g(t)$. 

Hence, \refe{continuous-K-eqn} and \refe{nfd-mobility-service} form a bathtub model of network trip flows served by mobility service vehicles. In this model, $\rho(t)$ is exogenously given. The corresponding model in the $(z,x)$ coordinates can be obtained.
Numerically, one can extend the simple difference and integration methods in \refe{intral-numerical-method} to solve $\la(t)$, $z(t)$, and $v(t)$ from the integral bathtub model, \refe{def:zt}, \refe{nfd-mobility-service}, and  (\ref{integral-lambda}a); then $K(t,x)$ can be solved from \refe{integral-bathtub}.

 Vickrey's bathtub model can be extended for mobility service systems when the trip distance follows a time-independent negative exponential distribution with $\Phi(t,x)=\tilde \Phi(t,x)=e^{-\frac xB}$:
 \bsq
\bqn
\dot \la(t)&=&f(t) -g(t),\\
g(t)&=&\frac 1 B\la(t) V(\rho(t),\la(t),f(t),g(t)).
\eqn
\esq
Then from \refe{nfd-mobility-service} and \refe{def:zt} we can calculate $v(t)$ and $z(t)$.
Note that Theorem \ref{thm:vickreymodel} is valid for any $v(t)>0$. The integral form of the extended Vickrey's bathtub model is still given by \refe{integral-Vickrey}, and \refe{G-Vickrey} can be used to calculate $G(t)$.

With deterministic trip distances, $\tilde B(t)$, \refe{continuous-K-eqn} can be simplified as
\bqn
\pd {}t K(t,x) -v(t) \pd{}x K(t,x)&=&f(t) H(\tilde B(t)-x),
\eqn
whose integral form is still \refe{integral-bathtub-deterministic}. The integral form of the conservation equation in $\la(t)$ is \refe{integral-la-deterministic}.

\subsection{Multi-commodity trips}

In \refe{K-model}, all vehicles travel at the same speed at the same time and make the same contribution to traffic congestion. But it can be extended for multi-commodity network trip flows, in which trips are grouped into commodities based on the characteristics of the trips and the serving vehicles:  different  commodities of trips can travel at different speeds, served by different types of vehicles, have different contributions to traffic congestion,  have different distances, or share with different numbers of trips; and different types of vehicles can have different contributions to traffic congestion (such as trucks and lane-changing vehicles), travel at different speeds, and have different occupancies (such as high-occupancy vehicles and buses).

For commodity $m$ ($m=1,\cdots, M$), we denote the number of active trips by $\la_m(t)$, the number of active trips with a remaining distance not smaller than $x$ by $K_m(t,x)$, the travel speed by $v_m(t)$, the cumulative travel distance by $z_m(t)$, the in-flux by $f_m(t)$, the proportion of the entering trips with distances not smaller than $x$ by $\tilde \Phi_m(t,x)$, and the probability density function of the active trips' remaining distances at $t$ by $\varphi_m(t,x)$. Then the generalized bathtub model can be written as
\bsq
\bqn
\dot z_m(t)&=&v_m(t),\\
\pd {}t K_m(t,x) -v_m(t) \pd {}x K_m(t,x) &=& f_m(t) \tilde \Phi_m(t,x),\\
\dot \la_m(t)&=&f_m(t)-\varphi_m(t,0) \la_m(t) v_m(t).
\eqn
The corresponding integral forms of $K_m(t,x)$ and $\la_m(t)$ are
\bqn
K_m(t,x)&=&K_m(0,x+z_m(t))+\int_0^t f_m(s) \tilde \Phi_m(s,x+z_m(t)-z_m(s))ds,\\
\la_m(t,x)&=&K_m(0,z_m(t))+\int_0^t f_m(s) \tilde \Phi_m(s,z_m(t)-z_m(s))ds.
\eqn
\esq
\refe{intral-numerical-method} can be applied to solve these models under general initial and boundary conditions. 
The model can be simplified in special cases with time-independent negative exponential and deterministic distributions of trip distances, and Theorem \ref{thm:vickreymodel} is valid for any commodity.\footnote{For multi-commodity trip flows, since different commodities can have different travel speeds and cumulative travel distances, the $(z_m,x)$ coordinates may not be synchronized. Thus the corresponding models in such coordinates may not be useful. }

The corresponding fundamental diagram of multi-commodity network trip flows can be written as
\bqn
v_m(t)&=&V_m(\vec \rho(t), \vec \la(t), \vec f(t), \vec g(t) ),
\eqn
where $\vec \rho(t)$ is the vector of different types of vehicles,  $\vec \la(t)$ the vector of the numbers of all commodities' active trips, $\vec f(t)$ the vector of the in-fluxes of all commodities, and $\vec g(t)$ the vector of the out-fluxes of all commodities.  The fundamental diagram may be different, when different types of cars have different passenger car equivalent values; the speed-density relations for different modes, such as metro and privately operated vehicles, could be uncorrelated; when different vehicles have different occupancies; when there are high-occupancy vehicle lanes or high-occupancy toll lanes; when there are freight trucks; when there are mobility service vehicles and privately operated vehicles.
Such fundamental diagrams of multi-commodity network trip flows have to be derived, calibrated and validated for different multi-commodity trip flow systems.

\section{Conclusion}
In this study we first presented a generalized bathtub model of network trip flow served by privately operated vehicles, in which a road network is treated as a single bathtub, and all vehicles' speeds at a time instant are identical and determined by the network fundamental diagram. After defining the number of active trips at $t$, $\la(t)$, that with a remaining distance not smaller than $x$, $K(t,x)$, the cumulative travel distance, $z(t)$, and other variables, we derived four equivalent differential formulations of the model from the conservation of trips. Then we defined and discussed the properties of stationary and gridlock states, derived the integral form of the bathtub model with the characteristic method, obtained equivalent formulations in the $(z,x)$ coordinates, and presented two numerical methods and an example. We further studied equivalent formulations and solutions for two special types of distributions of trip distances: negative exponential or deterministic. In particular, when the trip distance follows the time-independent negative exponential distribution, the model becomes Vickrey's bathtub model; and other equivalent formulations were obtained for constant trip distances. Finally we extended the fundamental diagram and the bathtub model for multi-commodity trip flows with trips served by mobility service vehicles.

The bathtub model captures the impacts on traffic congestion (in terms of $\la(t)$) of the demand, represented by the in-flux and entering trips' distances, and the supply, represented by the network fundamental diagram. For trips served by privately operated vehicles, the bathtub model is similar to the LWR model as it is also derived from the fundamental diagram and the conservation law of vehicles; however as a transport equation with a nonlocal speed, the bathtub model admits no shock or rarefaction waves, since the characteristic curves are much simpler without interactions among each other. Thus it is much simpler than the network kinematic wave models \citep{jin2012kinematic}, which requires detailed network topology as well as origin-destination and route demand data. The bathtub model is similar to the point queue model for a single bottleneck \citep{vickrey1969congestion}, which captures the impacts on the queue size of the demand represented by the in-flux and the supply by the capacity; however, the supply in the bathtub model decreases with more vehicles in the network, but that in the point queue model is usually assumed to be relatively constant even with capacity drop.
Note that the bathtub model is fundamentally different from the network fundamental diagram or the extended fundamental diagram of network trip flows, since the latter describes the stationary relationship between the speed and the number of active trips as well as the number of running vehicles, but the former combines the latter and the conservation of trips to describe the dynamics of trip flows. That is, the relationship between the bathtub model and the network fundamental diagram is similar to that between the LWR model and the link fundamental diagram. In this sense, \citep{godfrey1969mechanism} is a pioneering contribution to network traffic characteristics as \citep{greenshields1935capacity} is to link traffic characteristics, and \citep{vickrey1991congestion} is a pioneering contribution to network trip flow dynamics as \citep{lighthill1955lwr,richards1956lwr} is to link traffic flow dynamics.

The major contributions of this study are in the following. First, this study presented a unified framework for modeling network trip flows with general distributions of trip distances, including negative exponential, constant, and regularly sorting trip distances studied in the literature. Second, the integral form of the bathtub model was derived with the characteristic method and used to analytically solve the model under special initial and boundary conditions. Third, the bathtub model in the $z$ coordinate was derived to simplify the analyses. Fourth, two equivalent numerical methods were presented to solve the bathtub model: one based on the differential form, and the other based on the integral form. Fifth, six equivalent conditions for Vickrey's bathtub model to be applicable were derived, and it was clarified that Vickrey's bathtub model does not apply when all trips' distances are constant. Finally, the bathtub model was shown to apply when vehicles and trips are disintegrated and when there are different types of trips.

The bathtub model can be useful for network-level studies as it is based on simplified representations of the demand and supply. However, this study is only a starting point, and many further extensions are possible and warranted. 
\bi
\item We will be interested in comparing the bathtub model with more detailed network flow models, such as the network kinematic wave models. In particular we will be interested in quantifying the differences between these models. 

\item The bathtub model is probabilistic in nature. In this study we are only concerned with the expected values of $K(t,x)$. In the future we will be interested in its variations, which can be useful for understanding the reliability of a transportation network. 

\item In this study, the fundamental diagram of network trip flows is deterministic. In the future, we will be interested in incorporating accidents and stochasticity into the fundamental diagram. In addition, we will also be interested in developing microscopic version of the bathtub model, which can be used to obtain the statistical solutions through Monte Carlo simulations. 

\item We will be interested in empirically studying the trip distance distribution, the fundamental diagram of network trip flows, and the trip flow dynamics, which will be important for determining the applicability of model.

\item In this study, the bathtub model is developed for a single transportation network. However, a transportation system with different facilities (such as roads in different directions, high-occupancy-vehicle lanes, high-occupancy-toll lanes, parking lots) or different modes (such as bikes, cars, buses and subway) may be better treated as multiple interacting bathtubs, and the bathtub model in this study may need to be extended or modified.

\item Such bathtub models of network trip flow can be used to the impacts on congestion of travelers' choices in trip chaining, ride sharing, departure time, mode, parking, and roads as well as connected and autonomous vehicles. 

\item Furthermore, we will be interested in studying congestion pricing and other traffic management schemes based on the bathtub model. 

\item In addition, the model can be used to study the environmental impacts of transportation networks and the corresponding management and control strategies.
\ei

\section*{Acknowledgments}
I'd like to thank Wenbin Jin, Irene Martinez Josemaria, Koti Reddy Allu, Professors Pete Fielding and Michael Hyland for their stimulating and helpful discussions. The views and results are the author's alone.

\pdfbookmark[1]{References}{references}

\end {document}